\documentclass[reqno]{amsart}
\usepackage{hyperref}
\usepackage{geometry}
\usepackage{enumitem}
\usepackage{amsmath}
\usepackage{amssymb}
\usepackage{amsthm}
\usepackage{amsrefs}
\usepackage{amsfonts}
\usepackage{xcolor}
\usepackage{mathtools}
\usepackage{stackengine}

\makeatletter
\@namedef{subjclassname@2020}{%
  \textup{2020} MSC}
\makeatother

\hypersetup{
	colorlinks   = true, 
	urlcolor     = blue, 
	linkcolor    = purple, 
	citecolor   = blue 
}

\def\XXint#1#2#3{{\setbox0=\hbox{$#1{#2#3}{\int}$ }
\vcenter{\hbox{$#2#3$ }}\kern-.6\wd0}}

\makeatletter
\newcommand*{\rom}[1]{\expandafter\@slowromancap\romannumeral #1@}

\newcommand{\restr}{%
  \,\raisebox{-.260ex}{\reflectbox{\rotatebox[origin=br]{-90}{$\lnot$}}}\,%
}
\newcommand{\ind}{\protect\raisebox{2pt}{$\chi$}}

\newcommand*\Laplace{\mathop{}\!\mathbin\bigtriangleup}

\newcommand{\di}{\text{d}}
\newcommand{\diam}{\text{diam}}
\newcommand{\dist}{\text{dist}}

\newcommand{\R}{\mathbb{R}}

\newcommand{\Z}{\mathbb{Z}}
\newcommand{\Sphere}{\mathbb{S}}
\newcommand{\h}{\mathbb{H}}

\newcommand{\N}{\mathbb{N}}
\newcommand{\p}{\mathbb{P}}

\newcommand{\inv}{^{-1}}
\newcommand{\bthm}{\begin{thm}}
\newcommand{\ethm}{\end{thm}}
\newcommand{\bproof}{\begin{proof}}
\newcommand{\eproof}{\end{proof}}
\newcommand{\blem}{\begin{lem}}
\newcommand{\elem}{\end{lem}}
\newcommand{\brem}{\begin{rem}}
\newcommand{\erem}{\end{rem}}
\newcommand{\eeqn}{\end{equation}}
\newcommand{\eeqnn}{\end{equation*}}
\newcommand{\beqn}{\begin{equation}}
\newcommand{\beqnn}{\begin{equation*}}
\newcommand{\eprop}{\end{prop}}
\newcommand{\eexm}{\end{exm}}
\newcommand{\enexm}{\end{nexm}}
\newcommand{\ecor}{\end{cor}}
\newcommand{\bcor}{\begin{cor}}
\newcommand{\bexm}{\begin{exm}}
\newcommand{\bnexm}{\begin{nexm}}
\newcommand{\bprop}{\begin{prop}}
\newcommand{\bdefn}{\begin{defn}}
\newcommand{\edefn}{\end{defn}}
\newcommand{\benum}{\begin{enumerate}}
\newcommand{\eenum}{\end{enumerate}}
\newcommand{\Bsplit}{\begin{split}}
\newcommand{\Esplit}{\end{split}}

\title[Harmonic measures on the Bowditch boundary]{Harmonic measures on the Bowditch boundaries of groups hyperbolic relative to virtually nilpotent subgroups}

\begin{document}
\theoremstyle{plain}
\newtheorem{thm}{Theorem}[section]
\newtheorem{lem}[thm]{Lemma}
\newtheorem{prop}[thm]{Proposition}
\newtheorem{cor}[thm]{Corollary}

\theoremstyle{definition}
\newtheorem{defn}[thm]{Definition}
\newtheorem{exm}[thm]{Example}
\newtheorem{nexm}[thm]{Non Example}
\newtheorem{prob}[thm]{Problem}

\theoremstyle{remark}
\newtheorem{rem}[thm]{Remark}

\author{Debanjan Nandi}
\address{\textbf{Debanjan Nandi} \\
Faculty of mathematics and computer science, Weizmann Institute of Science, 234 Herzl Street, Rehovot 76100, Israel}
\email{debanjan.nandi@weizmann.ac.il}
\date{}
\thanks{The author gratefully acknowledges support from ISF Grant 1149/18.}

\begin{abstract} 

For a group hyperbolic relative to virtually nilpotent subgroups, on a cusped graph associated to the group, we construct a random walk whose Martin boundary is the Bowditch boundary of the group. Moreover, the harmonic measure is a conformal density corresponding to a hyperbolic Green metric and is exact dimensional on the Bowditch boundary. The latter equipped with a visual distance induced by the Green metric is an Ahlfors-regular metric measure space. The dimension is given in terms of the drift, a Green drift and the asymptotic entropy. The Patterson-Sullivan density for the action on the cusped graph in this case is doubling, its dimension is obtained by looking at cusp excursions of geodesics.

\end{abstract}

\maketitle


\section*{Introduction}
The free group in two generators has a representation in PSL$(2,\R)$, say $\Gamma_1$, with quotient for its action on $\h^2$, an infinite volume surface of constant sectional curvature $-1$ and genus one, and a different representation, say $\Gamma_2$, as a lattice with a finite volume quotient for its action on $\h^2$, a surface of genus one with a cusp. The limit set of the action on $\h^2$, in the first case is a Cantor set quasisymmetric to the visual boundary of the Cayley graph of the free group, and in the second case $\Sphere^1$.

The Martin boundary of a simple random walk on the Cayley graph of the free group, identified with the orbit of $\Gamma_1$ by a quasiisometry, is homeomorphic to the visual boundary, which is `singular' with respect to the Bowditch boundary of $\Gamma_2$ (homeomorphic to $\Sphere^1$, the limit set of its action on $\h^2$ as well as the Martin boundary of the hyperbolic Brownian motion, see \S \ref{prelim} for definitions). In seeking parallels with the relatively hyperbolic case $\Gamma_2$ above, that is, with geometrically finite actions of Kleinian groups with non-trivial parabolic subgroups more generally, our goal in this paper is the following: given a relatively hyperbolic group $\Gamma$, we want 
\benum 
\item to realize its Bowditch boundary as the Martin boundary of a random walk on some hyperbolic space associated to $\Gamma$; admitting a proper, isometric geometrically finite $\Gamma$ action, such that,
\item the harmonic measure of the walk is a (quasi)conformal density for some hyperbolic metric in the space.
\eenum
Roughly, with the example discussed above in mind, the idea is to construct a process but using the relatively hyperbolic structure of $\Gamma$ alone, that would mimic some aspects of a Brownian motion (on the convex hull of the limit set) in an ambient space. 
The space we avail is the (slightly modified) Groves-Manning cusped graph \cite{GM} associated to the group; and for the process, we consider a random walk, whose `increments' are Markov and have a bounded range.

Relatively hyperbolic groups as groups acting geometrically finitely on hyperbolic metric spaces, are generalizations of Kleinian groups with the metric property of acting geometrically finitely on the upper half space (see \S \ref{prelim} for the precise definitions). A weak cusped space $X$ for a relatively hyperbolic group $\Gamma$ (hyperbolic with respect to a finite collection of subgroups $\mathcal{P}$, often denoted as the pair $(\Gamma,\mathcal{P})$), is a proper $\delta$-hyperbolic metric space which admits a cusp-uniform action by $(\Gamma,\mathcal{P})$ (see Bowditch \cite{Bow}, Healy-Hruska \cite{HH}). The Bowditch boundary is a metric invariant for weak cusped spaces that have a `constant horospherical distortion' property (\cite{Bow}, \cite{HH}, see \S \ref{prelim}), which means roughly speaking, the horospheres, based at parabolic limit points, in going towards infinity contract the same amount in all `directions'. Weak cusped spaces that have this property will be called cusped spaces. The Groves-Manning graph is a cusped space; if $\Gamma$ acts as a lattice in a rank one symmetric space, then the symmetric space is a cusped space for $\Gamma$ (\cite{HH}); any two cusped-spaces of a relatively hyperbolic group $\Gamma$, are $\Gamma$-equivariantly quasiisometric. As the visual boundary of cusped spaces for $\Gamma$, the Bowditch boundary of $\Gamma$ has a canonical conformal gauge.

\subsection*{Main results} 
We discuss the results in two parts.

\textbf{Part 1:} Given a hyperbolic graph $X$, let $\mathcal{D}(X)$ be the class of quasiruled metrics (quasigeodesic metrics where triples of points taken from suitable quasigeodesics satisfy a coarse reverse-triangle inequality, see \S \ref{prelim}) on $X$ which are quasiisometric to the graph metric. We will define a cusped graph $X_\Gamma$ in \S \ref{construction}, associated to a relatively hyperbolic group $\Gamma$. Following is the basic result. 

\bthm\label{char}
Let $\Gamma$ be a finitely-generated non-elementary relatively hyperbolic group with finitely-generated parabolic subgroups $\mathcal{P}=\{H_i\}_{i=1}^l$. Let $X=X_\Gamma$ be its cusped graph. Then the following are equivalent:
\benum
\item there exists a random walk on $X$ with Martin boundary as the Bowditch boundary of $\Gamma$ such that the harmonic measure is an Ahlfors-regular, $\Gamma$-equivariant conformal density associated to a metric in $\mathcal{D}(X)$;
\item $H_i$ for $i\in\{1,\ldots,l\}$, are virtually nilpotent.
\eenum
\ethm

$(1)\implies(2)$ uses the Bonk-Schramm embedding theorem \cite{BS}. See the proof of \cite[Theorem 10.2]{BS}, where by an argument based on the Assouad embedding theorem \cite{As}, it is shown that a proper hyperbolic metric space whose visual boundary is doubling has \textit{bounded geometry}; which is equivalent in the case the space admits a properly discontinuous, geometrically finite action by a group of isometries, to the parabolic subgroups of the group with respect to the action being virtually nilpotent, a result in Dahmani-Yaman \cite{DY}. 

For the proof of $(2)\implies (1)$, we construct a random walk $(X,P)$ according to the geometry of the Groves-Manning cusped graph which jumps orbits. We can not apply Ancona's proof \cite{Anc} directly to a random walk on $X$ as the degree of vertices is unbounded. So the usual Harnack inequality for positive superharmonic functions and the familiar version of Ancona's inequality for hyperbolic graphs of uniformly bounded degree are not true. One way of dealing with this situation is choosing a suitable volume for $X$, denoted $m$, which in our case even forms a reversible Markov chain $(X,P,m)$, to which Ancona's machinery does apply. It follows as in Blacher\'e-H\"{a}issinsky-Mathieu \cite{BHM}, by hyperbolicity and quasiisometry with the graph metric of a suitable Green metric defined for $(X,P,m)$, that the harmonic measure is a conformal density. The Ahlfors-regularity of the harmonic measure follows from a shadow lemma for the harmonic measure and the Green metric. The situation requires we estimate shadows of balls centered at non-orbit points also which is possible due to the geometry of the Groves-Manning graph and Sullivan's shadow lemma for shadows of orbit points. The Ahlfors-regularity of the harmonic measure for the visual metric induced by the Green metric is in contrast to the fact that the Patterson-Sullivan density of $\Gamma$ for its action on $X$, is in general only doubling with respect to the visual metric induced by the graph metric.

\textbf{Part 2:} In the case when a hyperbolic space has symmetry, it may impose some regularity on the harmonic measure of suitable random walks. We study this interplay in our setting next. 

As a first step in this direction, for cusped graphs of groups hyperbolic relative to virtually nilpotent parabolics, we prove that drifts/rates of escape (for the graph and Green metrics) and asymptotic entropy exist for the random walk constructed. This is done by looking at a `quotient walk' which we think of as the induced random walk on the quotient space under the action of $\Gamma$ on $X$, which can be coded symbolically as a process with Markov increments. The corresponding left shift has an invariant Markov probability measure, which makes it even a mixing transformation. The Kingman subadditive ergodic theorem applies. Write $\rho_X$ for the graph metric on the cusped graph $X$. For a reversible Markov chain $(X,P,m)$ on $X$, write $\p_m=\sum_{x\in X}m(x)\cdot \p_x$ ($m$ being the volume in $X$), for a measure in $X^{\N}$. Write $\rho_G$ for a suitable Green metric (defined in \S \ref{green_metric}) in $X$. 
\bthm\label{drift} Let $X$ be the cusped graph for a group, hyperbolic relative to virtually parabolic subgroups. There is a reversible Markov chain $(X,P,m)$, such that for $\{Y_i\}_{i\in\N_0}$, a process corresponding to $(X,P,m)$ the following limits exist and are constants $\p_m$-a.e $\omega\in X^\N$:
$$\hspace{0.4cm}(\text{drift})\hspace{3cm} \underset{n\to\infty}{\lim}  \frac{\rho_X(Y_0(\omega),Y_n(\omega))}{n}=:l,$$
$$(\text{Green drift})\hspace{2.5cm} \underset{n\to\infty}{\lim}  \frac{\rho_G(Y_0(\omega),Y_n(\omega))}{n}=: l_G,$$
$$(\text{asymptotic entropy})\quad\underset{n\to\infty}{\lim} \frac{-\log \,p^{(n)}(Y_0(\omega),Y_n(\omega))}{n}=:h.$$ 
\ethm

Next, we show that the harmonic measure for $(X,P,m)$ is exact dimensional, that is, the ratio of the logarithm of the harmonic measure of a ball and the logarithm of the radius of a ball (in suitable visual metrics) has a limit as the radius approaches zero (a detailed discussion of the literature on this topic in related contexts appears in Tanaka \cite{Tan}). The formula for the dimension can be given in terms of drift and asymptotic entropy. There are two approaches to this computation. One goes via a shadow lemma for the Green metric with respect to the harmonic measure using the Ancona inequality. This, and a sublinear geodesic tracking property for the random walk with respect to the Green and graph metrics leads to a formula for the dimension. This approach is parallel to \cite{BHM}. The main issue in our case is to extend the shadow lemma to the shadow of balls centered at non-orbit points. Another approach from Kaimanovich \cite{Kai1}, Le Prince \cite{LeP}, Tanaka \cite{Tan}, applicable rather generally for groups acting with exponential growth on proper hyperbolic spaces and jumps given by measures with finite first moment (see also Dussaule-Yang \cite{DuYa} for harmonic measures of random walks on Cayley graphs of relatively hyperbolic groups), utilizes properties of the asymptotic entropy, drift (Kaimanovich \cite{Kai2}; Maher-Tiozzo \cite{MT} in the non-proper case) in the hyperbolic setting and the ergodicity of the action of $\Gamma$ with respect to the harmonic measure in the boundary (for the lower bound on the dimension, \cite{Tan}). The dimension of $\mu$, the $D_\Gamma$-dimensional Patterson-Sullivan density of $\Gamma$ for its action on $X$, where $D_\Gamma$ is the critical exponent of $\Gamma$ acting on $X$ ($\mu$ being unique up to constant multiples as it is non-atomic, doubling and ergodic for the $\Gamma$-action, see Remark \ref{ergodicity} and Remark \ref{imp_2} below) is computed by an approach of Stratmann-Velani \cite{SV}, using a shadow lemma to study cusp excursions into horoballs. 

For a random walk $(X,P)$ let $\nu$ denote the harmonic measure. We denote by $B_{x,\infty}^X(\xi,r)$ and $B_{x,\infty}^G(\xi,r)$ balls in $\partial X$, centered at $\xi\in\partial X$ of radius $r>0$, respectively for visual metrics (with base-point $x\in X$, suppressed in notation) $d_{X,\epsilon_X}$ and $d_{G,\epsilon_G}$, induced by hyperbolic metrics $\rho_X$ (graph) and $\rho_G$ (Green) of $X$, with $\epsilon_X, \epsilon_G>0$ being suitable parameters for the visual metrics (see \S \ref{prelim} for definitions). Let $\mu$ be a Patterson-Sullivan density of dimension $D_\Gamma$, which is obtained by the Patterson construction.
\bthm[Drift, entropy and dimension]\label{dimdriftent}
Let $X$ be the cusped graph for a group, hyperbolic relative to virtually parabolic subgroups. There is a random walk $(X,P)$ so that for $x\in X$, we have for the harmonic measure $\nu_x$ supported in $\partial X$, $$\lim_{r\to 0}\frac{\log (\nu_x(B_\infty^X(\xi,r)))}{\log r}=\frac{l_G}{\epsilon_X\cdot l},\quad \lim_{r\to 0}\frac{\log (\nu_x(B_\infty^G(\xi,r)))}{\log r}=\frac{1}{\epsilon_G},$$ and
$$\lim_{r\to 0}\frac{\log (\nu_x(B_\infty^X(\xi,r)))}{\log r}=\frac{h}{\epsilon_X\cdot l},\quad \lim_{r\to 0}\frac{\log (\nu_x(B_\infty^G(\xi,r)))}{\log r}=\frac{h}{\epsilon_G\cdot l_G}.$$ For the Patterson-Sullivan density $\mu$ of $(\Gamma,X)$, we have $$\lim_{r\to 0}\frac{\log (\mu_x(B_\infty^X(\xi,r)))}{\log r}=\frac{D_\Gamma}{\epsilon_X}.$$
\ethm

See \S \ref{guiver} for a discussion on the relation between entropy and drift. The inequality  $h\leq l\cdot D_\Gamma$ (Guivarc'h \cite{Gui1}, Kaimanovich \cite{Kai3}, Ledrappier \cite{Led}, Vershik \cite{Ver}), is true for the walk we consider. Moreover a characterization of the equality as in \cite[Theorem 1.5]{BHM} also holds. We note a formula for the Green drift (Proposition \ref{intrep}) in terms of the Busemann-Martin cocycle (see Gou\"ezel-Math\'eus-Maucourant \cite{GMM} for the case of random walks with independent increments on hyperbolic groups).

\subsection*{Short survey}
Relatively hyperbolic groups were introduced in Gromov \cite{Gr2}. Further foundational work came in Farb \cite{Fa}, Bowditch \cite{Bow} where alternative definitions were given. Alternative characterizations followed in Osin \cite{Os}, Yaman \cite{Yam} and Groves-Manning \cite{GM}. In the recent work Healy-Hruska \cite{HH}, a characterization is given for the quasiisometric class of the Groves-Manning cusped graph. Relatively hyperbolic groups with virtually nilpotent parabolics were characterized in Dahmani-Yaman \cite{DY} as groups which have geometrically finite actions on hyperbolic spaces with bounded geometry, and appear in the recent work Mackay-Sisto \cite{MS} where a Bonk-Schramm-type characterisation \cite{BS} for them is proved in terms of embeddings into truncated real hyperbolic spaces.

Fundamental works on the topic of random walks in general hyperbolic groups and graphs are in Ancona \cite{Anc}, for the Martin boundary, and in Kaimanovich \cite{Kai2}, for the Poisson boundary. Ancona's inequality proved and used in \cite{Anc}, originally for transition probabilities with finite support is crucial for us.
Our paper draws motivation from Blacher\`{e}-H\"{a}issinsky-Mathieu \cite{BHM} where the authors study the properties of a Green metric, introduced in Blacher\`e-Brofferio \cite{BB}, proving its hyperbolicity and its consequences. 

There has recently been significant advancement in the understanding of random walks on relatively hyperbolic groups; see Dussaule \cite{Dus}, Dussaule-Gekhtman-Gerasimov-Potyagailo \cite{DGGP}, Gekhtman-Gerasimov-Potyagailo-Yang \cite{GGPY}. There are now, newer approaches to the Ancona inequality, see Gouezel-Lalley \cite{GL}, Gouezel \cite{Gou1}, \cite{GGPY}, that work for infinite range jumps of finite super-exponential moment in Cayley graphs of hyperbolic and relatively hyperbolic groups.

The question of ascertaining the dimension of the harmonic measure in related settings has been treated in Ledrappier \cite{Led1, Led2}, Kaimanovich \cite{Kai1}, Le Prince \cite{LeP}, \cite{BHM}, Tanaka \cite{Tan}. It has been treated for random walks on relatively hyperbolic groups in \cite{DuYa}, for harmonic measures in both Floyd and Bowditch boundaries. The random walks considered in the above works are walks that jump along orbits. The approach to dimension in \cite{Tan} applies also to our setting with appropriate modification as mentioned earlier-the upper bounds and lower bounds are computed separately; the upper bound using a method in \cite{LeP}, the lower bound is obtained using ergodicity of the $\Gamma$-action with respect to the harmonic measure on $\partial X$. Another approach to dimension for the harmonic measure available to us, runs parallel to the one in \cite{BHM} of using the hyperbolicity of a Green metric and of Stratmann-Velani \cite{SV}, for the Patterson-Sullivan measure which uses excursions of geodesics for geometrically finite actions. 

Random walks along orbits behave differently than the walk considered here which also jump orbits. The properties of the Bowditch boundary studied here are not captured using random walks pushed-forward from a Cayley graph. In a recently announced paper, the results in \cite{SV} on cusp-excursions and a general shadow lemma (in \cite{SV} for geometrically finite Kleinian group actions and Schapira \cite{Sch} in the case of geometrically finite actions on simply connected manifolds of negative curvature) have been generalized to the setting of certain projective spaces in Bray-Tiozzo \cite{BT}; one application of which is to show in fact that suitable conformal densities corresponding to groups acting geometrically finitely are singular with respect to the harmonic measures corresponding to random walks along orbits. See \cite{GT} for results comparing a class of quasiconformal measures and harmonic measures corresponding to random walks along orbits (proving singularity) in the CAT($-1$) setting. 

\subsection*{Random walks in graphs with discrete group actions} 
Consider a graph $X$ with action by an infinite discrete group $\Gamma$ of isometries (automorphisms). Given a probability measure $\mu$ on $\Gamma$, there is a random walk on $X$ driven by $\mu$, that is a random walk for each $x\in X$, starting at $x$, $Z_n\cdot x=w_n\cdots w_1\cdot x$, where $w_i$ are $\Gamma$-valued independent random variables with distribution $\mu$ (whose support is assumed to generate $\Gamma$ as a semigroup). These well studied (for example \cite{Kai2}, \cite{MT}) random walks come within the purview of $\Gamma$-invariant Markov chains with the larger state space $X$ (see \S \ref{gen_walk}). Allowing for random walks that jump orbits are useful from a geometric or group-theoretic point of view as they give useful information about the space and the group.

\subsection*{Structure of the paper} In \S \ref{prelim} we collect the relevant definitions and present some preliminary facts and state a lemma needed in identifying the Martin boundary of our walk. In \S \ref{construction} we construct a random walk on a cusped graph for a relatively hyperbolic group with virtually nilpotent parabolic subgroups and prove an isoperimetric inequality for it. In \S \ref{BowditchisMartin} we show that the Martin boundary of the random walk constructed in \S \ref{construction} is the Bowditch boundary of the group. In \S \ref{metricmeasurestructure} we introduce a Green metric for the random walk and show that the harmonic measure of the random walk is a conformal density with respect to the Green metric. We show using a shadow lemma for the harmonic measure and the Green metric that the harmonic measure is Ahlfors-Regular (for the visual metric induced by the Green metric), and obtain the growth of orbits with respect to the Green metric. In \S \ref{dimensionharmonic} we look at a recurrent quotient walk to obtain drifts/rates of escape (for the graph and Green metrics) and an asymptotic entropy for the random walk. We use cusp excursions of geodesics with respect to the Patterson-Sullivan measure to obtain its dimension. We use a geodesic tracking property of the walk to obtain the dimension of the harmonic measure and relate it to the drift and asymptotic entropy. 

We write $A\lesssim B$ to mean $A\leq C\cdot B$, for some $C>0$ which can be computed independently of $A,B$. Similarly for $A\gtrsim B$. By $A\approx B$ we mean simultaneously $A\lesssim B$ and $B\lesssim A$. 

\subsection*{Acknowledgements} The author is grateful to Omri Sarig for helpful discussions and support. The author is grateful to Shubhabrata Das and Suraj Krishna M. S. for a question that provided the initial motivation for this research and for helpful discussions. The author thanks the International Centre for Theoretical Sciences, Bengaluru as part of the programme \emph{Probabilistic Methods in Negative Curvature} which acquainted him with this area.

\section{Preliminaries}\label{prelim}

\subsection{Doubling metric measure spaces}
A metric measure space is a metric space $(A,d,\mu)$ with a Borel regular measure $\mu$. For a constant $C\geq 1$, $(A,d)$ is $C$-doubling if every ball of radius $0<r<\diam(A)$ can be covered by at most $C$ balls of radius $r/2$. For $C\geq 1$ a metric measure space $(A,d,\mu)$ is $C$-doubling if for every $0<r<\diam(A)$, $\xi\in A$, $$\mu(B(\xi,2\cdot r))\leq C\cdot \mu(B(\xi,r)),$$ where $B(\xi,r)$ is a ball of radius $r>0$ in $(A,d)$, centered at $\xi$.

Next we define a class of maps of metric spaces, that distort distances with uniform bounds, at each scale (Tukia-Vais\"al\"a \cite{TV}, Heinonen \cite{Hei}).
\bdefn[Quasisymmetry] A homeomorphism $f:(A_1,d_1)\longrightarrow(A_2,d_2)$ is a $\varphi$-quasisymmetry given an increasing homeomorphism $\varphi:(0,\infty)\to (0,\infty)$, if for $\xi,\eta,\zeta\in A_1$ and $t>0$,
$$d_1(\xi,\eta)\leq t\cdot d_1(\eta,\zeta)\implies d_2(f(\xi),f(\eta))\leq \varphi(t)\cdot d_2(f(\eta),f(\zeta)).$$
\edefn
\brem
The inverse of a quasisymmetry is a quasisymmetry.
\erem
\blem[Quasisymmetry and doubling] If for $C>0$, $(A_1,d_1)$ is a $C$-doubling metric space and $f:(A_1,d_1)\longrightarrow(A_2,d_2)$ is a $\varphi$-quasisymmetry, then $(A_2,d_2)$ is $C'=C'(C,\varphi)$ doubling. 

If $(A_1,d_1,\mu_1)$ is $C$-doubling metric measure space, and $f:(A_1,d_1)\longrightarrow(A_2,d_2)$ is a $\varphi$-quasiisometry, then $(A_2,d_2, \mu_1)$ is $C'=C'(C,\varphi)$ doubling. 
\elem

\subsection{Hyperbolic spaces} 
We define some of the terms needed in the paper, see Bridson-Haefliger \cite{BH} for details on the basic notions surrounding hyperbolic metric spaces.
\bdefn[Gromov product] For a metric space $(Y,\rho)$ and a point $y\in Y$, the Gromov product with respect to $y\in Y$ is the function $(\cdot|\cdot)_y: Y^2\to[0,\infty)$ defined by
$$(u|v)_y=\frac{1}{2}\cdot (\rho(y,u)+\rho(y,v)-\rho(u,v)),$$ for $u,v\in Y$.
\edefn

\bdefn[Hyperbolic metric space]
Let $\delta\geq 0$. A metric space $(Y,\rho)$ is $\delta$-hyperbolic if there exists $y\in Y$, such that for all $u, v, w\in Y$, 
$$(u|w)_y=\min\{(u|v)_y,(v|w)_y\}-\delta.$$
\edefn
 
A metric space $(Y,\rho)$ is called \textit{geodesic} if for any pair $x,y\in Y$ there is an isometry from an interval of $\R$ to $Y$ with the endpoints of the interval being mapped to $\{x,y\}$. 
As a convention we often write $[x,y]$ to denote a geodesic joining $x$ to $y$, and $[x,y,z]$ for a geodesic triangle with vertices $x,y,z\in Y$.
\brem In a $\delta$-hyperbolic geodesic metric space, it is convenient to think of the Gromov product $(u|v)_y$ coarsely as the distance from $y$ of any geodesic $[u,v]$ joining $u, v$, as the following is true: $|\dist_{\rho}(y,[u,v])-(u|v)_y|\leq \delta$.
\erem

(Graph-metric). 
For a graph $Y=(V(Y),E(Y))$ with vertex set $V(Y)$ and edge set $E(Y)$, consider the $1$-simplex also denoted $Y$, such that every edge $e=(e_1,e_2)\in E$ for $e_1,e_2\in V$, is identified isometrically to the unit interval $[0,1]\subset\R$. This structure as a $1$-simplex induces a unique length metric on $Y$, such that every pair of points $x,y\in V(Y)$ has distance equal to the minimum number of edges required to go from $x$ to $y$, which we denote by $\rho_Y$ and call the graph metric of $Y$.

\brem
Note that $(Y,\rho_Y)$ is a geodesic metric space.
\erem

\bdefn[Quasiisometry]
For metric spaces $(Y,\rho_1)$, $(Y,\rho_2)$, and numbers $K\geq 1, L>0$, a $(K,L)$-quasiisometric embedding $f:Y_1\to Y_2$ is a map such that 
$$\frac{1}{K}\cdot \rho_1(u,v)-L\leq \rho_2(f(u),f(v))\leq K\cdot\rho_1(u,v)+L,$$ for all $u,v\in Y$. For $C\geq 0$, the map $f$ is a $(K,L,C)$-quasiisometry if $f$ is a $(K,L)$-quasiisometry and there is a map $g:Y_2\to Y_1$, such that $\rho_1(y,g\circ f(y))\leq C$, for all $y\in Y$.
\edefn

For a metric space $(Y,\rho)$ and $K\geq 1$, $L\geq 0$, a $(K,L)$-\textit{quasigeodesic} is a $(K,L)$-quasiisometric embedding from an interval of $\R$ to $X$.

\brem[Morse Lemma] Any $(K,L)$-quasigeodesic in a geodesic hyperbolic space lies in an $0<D=D(\delta,K,L)$-neighbourhood of a geodesic.\erem

The following notion from \cite[p. 8]{BHM}, is needed.
\bdefn[Quasiruling]
For a metric space $(Y,\rho)$, and numbers $K\geq 1, L\geq 0, \tau\geq 0$, a $(\tau,K,L)$-quasiruler is a $(K,L)$-quasigeodesic $\gamma:I\to X$ such that for all $s<t<u\in I$,
$$\rho(\gamma(s),\gamma(t))+\rho(\gamma(t),\gamma(u))-\rho(\gamma(s),\gamma(u))\leq 2\tau.$$ 
A metric space $(Y,\rho)$ is quasiruled if there exist $\tau,K,L$ such that every pair of points in $Y$ can be joined by a $(K,L)$-quasigeodesic and every $(K,L)$-quasigeodesic is a $\tau$-quasiruler.
\edefn

\brem[Thin triangles condition] For geodesic metric spaces there is an alternative criterion for hyperbolicity. Let $\delta\geq 0$. A geodesic triangle in a geodesic metric space $(Y,\rho)$ is $\delta$-thin, if every edge is contained in the $\delta$-neighbourhood of the other two. A geodesic metric space is $\delta$-hyperbolic if and only if every geodesic triangle is $\delta$-thin.
This criterion is useful, for instance, in application to the result below \cite[Theorem A.1]{BHM} which we use.
\erem

\bthm[Quasiisometry and hyperbolicity]\label{QI2H} Let $(Y_1,\rho_1)$ be a geodesic hyperbolic metric space, $(Y_2,\rho_2)$ a metric space and $\phi:Y_1\to Y_2$ be a quasiisometry. Then $(Y_2,\rho_2)$ is hyperbolic if and only if it is quasiruled (the constants can be bounded using data given).
\ethm

\bdefn[Visual boundary]\label{visual_bnd}
For a metric space $(Y,\rho)$, a sequence $\{z_n\}_n\subset Y$ converges at infinity if $(z_i|z_j)_y\xrightarrow{i,j\to\infty}\infty$, for some $y\in Y$. Sequences $\{z_n\}_n, \{z'_n\}_n\subset Y$ are defined equivalent if 
$(z_i|z'_j)_y\xrightarrow{i,j\to\infty}\infty$, for some $y\in Y$. The set of equivalence classes of sequences converging at infinity, with topology generated by sets of the form $$U(\{z_n\}_n,c)=\left\lbrace[\{z'_n\}_n]\,\middle|\,\{z'_n\}_n \;\text{convergent at infinity},\;\underset{i,j\to\infty}{\limsup}\,(z_i|z'_j)_y\geq c,\, y\in Y\right\rbrace,$$ for $c>0$ and $\{z_n\}_n$ convergent at infinity, is the visual boundary of $Y$.
\edefn

The Gromov-product can be extended to the boundary: for $\xi,\eta\in \partial X$, one defines $$(\xi|\eta)_y=\sup \,\left\lbrace \underset{i,j\to\infty}{\liminf} \,(z_i|z'_j)_y\,\middle|\,\{z_n\}_n\in\xi, \{z'_n\}_n\in\eta\right\rbrace.$$

\bthm[Visual metric]\label{visual_met} Let $(Y,\rho)$ be $\delta$-hyperbolic. There exists $\epsilon_0=\epsilon_0(\delta)$, such that for all $0<\epsilon<\epsilon_0$, there exists a visual metric $d_{y, \epsilon,\rho}$ on $\partial X$, such that $$\frac{1}{C}\cdot e^{-\epsilon\cdot(\xi|\eta)_y}\leq d_{y, \epsilon,\rho}(\xi,\eta)\leq C\cdot e^{-\epsilon\cdot(\xi|\eta)_y},$$ for some $C=C(\epsilon,\delta)$ and all $\xi,\eta\in Y$. If $Y$ is proper, $d_{y, \epsilon,\rho}$ generates a compact topology equivalent to the one in Definition \ref{visual_bnd}.\ethm
The dependence of $d_{y, \epsilon,\rho}$ on $y\in Y$ is bilipschitz and we will usually drop it from the notation for the visual metric when the point $y$ is fixed. We will denote $\{z_n\}_n\in\xi\in\partial Y$, by $z_n\to\xi$.
\bdefn[Busemann function] For a hyperbolic metric space $(Y,\rho)$,  the Busemann function $\beta:\partial X\times X^2\to \R$, is defined
$$\beta_\xi(u,v)=\sup\,\left\lbrace\underset{n\to\infty}{\limsup} \,(\rho(u,z_n)-\rho(v,z_n)) \,\middle|\,\{z_n\}_n\subset Y, \,z_n\to \xi\right\rbrace,$$
for $\xi\in\partial Y$ and $u,v\in Y$. 
\edefn

\bdefn[Shadows] For a $\delta$-hyperbolic space $(Y,\rho)$, $R>0$ $y_0, y\in Y$, the shadow of $B_{\rho}(y,R)$ cast from $y_0$, is the set $$\hat{S}_\rho(y_0, B_\rho(y,R))=\{\xi\in\partial Y\,\mid \, (\xi|y)_{y_0}\geq \rho(y_0,y)-R\}.$$
\edefn

A proper $\delta$-hyperbolic space $(Y,\rho)$ is called $R$-\textit{visual} if there exists $R>0$ such that for each $y,z\in Y$, there is a geodesic ray starting at $y$ that intersects $B(z,R)$. A proper $(K,L,\tau)-$quasiruled $\delta$-hyperbolic space $(Y,\rho)$ is $(K,L,\tau,R)$-\textit{quasivisual} if there is $R>0$, such that for each $y,z\in Y$, there is a $(K,L,\tau)$-quasiruled ray starting at $y$ that intersects $B(z,R)$.

For a quasivisual space $(Y,\rho)$, shadows can also be defined as follows: for  $R>0$ $y_0\in Y$, the shadow of $E\subset Y$ cast by $y_0$, is the set \[S_\rho(y_0, E)=\left\lbrace \xi\in\partial Y\,\middle| \,
\begin{tabular}{@{}l@{}}
$\exists\;\text{quasiruled ray starting at}\;y_0,\;\gamma_0,\;\text{such that}$ \\
 $\gamma_0(\infty)=\xi, \;\gamma_0\cap E\neq \emptyset$  \end{tabular}
 \right\rbrace.
 \]
Then there is the following result; see \cite[Proposition B.5]{BHM}. 
\bthm
Let $(Y,\rho)$ be a $(K,L,\tau,R)$-quasivisual space. There exist $R_0=R_0(K,L,\tau,R)>0, C=C(K,L,\tau,R)>0$, such that for all $y_0,y\in Y$, $r>R_0$, 
$$S_\rho(y_0, B_\rho(y,r-C))\subset \hat{S}_\rho(y_0, B_\rho(y,r))\subset S_\rho(y_0, B_\rho(y,r+C)).$$
\ethm
This leads to the following needed fact \cite[Proposition 2.1]{BHM}.
\blem
Let $(Y,\rho)$ be $(K,L,\tau,R)$-quasivisual space. Then for any $\tau'\geq 0$, there exist $R_0=R_0(K,L,\tau,R,\tau')>0, C=C(K,L,\tau,R,\tau')>0$, such that for all $r>R_0$, $\xi\in \partial Y$ and $y_0, y\in Y$, with $(y_0|\xi)_y\geq \tau$, it holds for a visual metric $d=d_\epsilon$ in $\partial Y$:
$$B_d\left(\xi, \frac{1}{C}\cdot e^{r\epsilon}\cdot e^{-\epsilon\rho(y_0,y)}\right)\subset S_\rho(y_0, B_\rho(y,))\subset B_d(\xi, C\cdot e^{r\epsilon}\cdot e^{-\epsilon\rho(y_0,y)}).$$
\elem

Next we mention the duality between quasiisometries of the space and quasisymmetries of the boundary. See Buyalo-Schroeder \cite[Theorem 5.2.17]{BuySch} for the following result. 
\bthm[QI to QS]\label{QItoQS} Given proper visual $\delta$-hyperbolic spaces $(Y_1,\rho_1)$ and $(Y_2,\rho_2)$, a quasiisometry $\phi: Y_1\longrightarrow Y_2$ induces a $\varphi$-quasisymmetry $\hat{\phi}:(\partial Y_1, d_{\epsilon_1,\rho_1})\longrightarrow (\partial Y_2,d_{\epsilon_2,\rho_2})$, between their visual boundaries, where $\varphi$ depends only on $\phi, \delta, \epsilon_1, \epsilon_2$.
\ethm
\brem[Conformal gauge]
For a proper $\delta$-hyperbolic metric space, the family of metrics quasisymmetric to the visual metric in the visual boundary (Theorem \ref{visual_met}), is called the \textit{conformal gauge} of the visual boundary. Under a natural assumption (uniform perfectness) for the metric, the converse to Theorem \ref{QItoQS} holds. In particular, every quasisymmetry of a visual metric can be extended to a quasiisometric self map of the hyperbolic space \cite[Corollary 7.2.3]{BuySch}.
\erem

\subsection{Relatively hyperbolic groups} 

If a discrete group $\Gamma$ acts by isometries on a $\delta$-hyperbolic space $Y$, there is the induced action on the visual boundary $\partial Y$: an isometric image of a sequence convergent at infinity is also convergent at infinity and gives a well-defined action on the visual boundary. We call the set of accumulation points of the action in $\partial Y$ the \textit{limit set} of $\Gamma$ and denote it $\Lambda_\Gamma(Y)$ (usually $\Lambda_\Gamma$ when $Y$ is fixed).
We provide a definition of conical limit points that will be useful for us. See \cite[Proposition A.2]{HH} for equivalence with another familiar definition of conical limit points for general convergence group actions on compact metrizable spaces.
\bdefn[Conical limit points] If $\Gamma$ acts by isometries on a proper geodesic $\delta$-hyperbolic space $(Y,\rho)$ such that $\partial Y$ has cardinality at least three, then a point $\xi\in\Lambda_\Gamma(Y)$ is conical if there exists a geodesic ray $\gamma$ with end point $\xi$, a point $y\in Y$, and a sequence $\{g_i\}_i\subset\Gamma$ such that $g_i(y)\to\infty$ and $\rho(g_i(y),\gamma)$ is uniformly bounded. 
\edefn
\bdefn[Parabolic limit points] If $\Gamma$ acts by isometries on a proper $\delta$-hyperbolic space $(Y,\rho)$, then a point $\xi\in\Lambda_\Gamma$ is parabolic, if its stabilizer $P_\xi$ in $\Gamma$ for the induced action of $\Gamma$ in $\partial Y$, has $\xi$ as a unique limit point in $\partial Y$. Any subgroup of $P_\xi$ is called parabolic, and $P_\xi$ is called a maximal parabolic subgroup of $\Gamma$. A parabolic limit point $\xi\in\partial Y$ is called bounded parabolic, if the action of $P_\xi$ on $\Lambda_\Gamma(Y)\setminus\{\xi\}$ is cocompact. 
\edefn

Let $\Gamma$ be a finitely generated group and $\mathcal{P}$ a finite collection of subgroups of $\Gamma$.
\bdefn[Cusp-uniform action] The pair $(\Gamma,\mathcal{P})$ has a cusp-uniform action on a proper $\delta$-hyperbolic metric space if the action is:
\benum
\item properly discontinuous by isometries
\item such that $\mathcal{P}$ is a set of representatives of conjugacy classes of maximal parabolic subgroups,
\item cocompact on the complement of a $\Gamma$-equivariant family of open horoballs.
\eenum
The space $X$ admitting a cusp-uniform action by a pair $(\Gamma,\mathcal{P})$ is called a weak-cusped space for the pair $(\Gamma,\mathcal{P})$ \cite{HH}.
\edefn

\bdefn[Relative hyperbolicity] We call a finitely generated group $\Gamma$ relatively hyperbolic if there is a finite collection of subgroups $\mathcal{P}$ and a proper, visual $\delta$-hyperbolic space which is a weak-cusped space for the pair $(\Gamma,\mathcal{P})$.
\edefn
\brem[Geometrical finiteness] A cusp-uniform action of a relatively hyperbolic pair $(\Gamma, \mathcal{P})$ on a weak-cusped space $Y$ is geometrically finite: $\partial Y$ constitutes of conical limit points and bounded parabolic limit points \cite{Bow}, \cite{HH}. 
\erem
Examples include free groups, fundamental groups of complete finite volume manifolds of pinched negative sectional curvature; these examples are also hyperbolic relative to virtually nilpotent subgroups. See Dahmani \cite{Dah} for general ways of constructing relatively hyperbolic groups.

Given a relatively hyperbolic group $(\Gamma,\mathcal{P})$, any two weak-cusped spaces have $\Gamma$-equivariantly homeomorphic visual boundaries \cite{Bow}.
\bdefn[Bowditch boundary]The Bowditch boundary of a relatively hyperbolic pair $(\Gamma,\mathcal{P})$ is the visual boundary of a weak-cusped space for $(\Gamma,\mathcal{P})$ (as a topological space).
\edefn
\brem[Cusped spaces and quasiisometry] There are the following remarks:
\benum
\item The Bowditch boundary is $\Gamma$-equivariantly homeomorphic for two weak cusped spaces $Y_1$ and $Y_2$ for $(\Gamma,\mathcal{P})$. There may not be a quasiisometry between $Y_1$ and $Y_2$ which is $\Gamma$-equivariant (Healy \cite{Hea}) and consequently, by uniform perfectness of visual metrics on boundaries of weak cusped spaces \cite{HH}, there exist visual metrics on the Bowditch boundary which are not $\Gamma$-equivariantly quasisymmetric \cite{HH}.
\item Cusped spaces are weak cusped spaces that have the following constant horospherical distortion property: for $H\in\mathcal{P}$, and any $y\in Y$, there is a bilipschitz function $\psi_{y,H}:(0,\infty)\to (0,\infty)$, such that for $h,h'\in H$, $e^{\psi(\rho(h(y),\rho(h'(y)))}\approx d_H(h,h')$, with an absolute constant, where $d_H$ is a Cayley graph metric for a fixed choice of finite generators. Cusped spaces $Y_1$ and $Y_2$ for a relatively hyperbolic group $(\Gamma,\mathcal{P})$ are $\Gamma$-equivariantly quasiisometric \cite{HH}.
\item In \S \ref{construction} we present (a minor variant of) the construction from \cite{GM} which is a cusped space. 
\eenum
\erem

We note the following theorem \cite{Bass}, \cite{Gui}, \cite{Gr1}.
\bthm[Bass, Guivarc'h, Gromov] Finitely generated groups of polynomial growth and finitely generated virtually nilpotent groups are the same class. The cardinality of a ball of radius $n$ is comparable by absolute constants to $n^{d_H}$, where $d_H$ is the rank of a finite index nilpotent subgroup of $H$.
\ethm
The growth function of a finitely generated group $H$ for some fixed set of generators will be denoted $g_H$.

\subsection{Quasiconformal densities} 
\bdefn[QC densities] For a proper $\delta$-hyperbolic space $Y$, a $\Gamma$-quasiconformal density $\nu=\{\nu_y\}_{y\in Y}$, of dimension $\alpha>0$ is a map $y\mapsto \nu_y$ from $Y$ to the set of finite Borel measures on $\partial Y$, such that the following two properties are satisfied:
\benum
\item ($\Gamma$-equivariance) For all $y\in Y$, $g\in \Gamma$, $g_\ast \nu_y= \nu_{gy}$, where $g_\ast\nu_y$ is the pushforward $\nu_y\circ g\inv$.
\item (Quasiconformality) There is $C\geq 1$, such that for all $x,y\in Y$, $\nu_y$-a.e $\xi\in \partial Y$, $$\frac{1}{C}\cdot e^{-\alpha\cdot\beta_\xi(x,y)}\leq \frac{d\nu_x}{d\nu_y}(\xi)\leq C\cdot e^{-\alpha\cdot\beta_\xi(x,y)}.$$
\eenum
If $C$ can be chosen one, $\nu$ is called a $\Gamma$-conformal density.
\edefn
Examples are Patterson-Sullivan measures on the boundaries of Cayley graphs of hyperbolic groups, of non-elementary discrete groups of isometries acting on hyperbolic spaces, harmonic measures of suitable random walks on Cayley graphs of hyperbolic groups, and as we shall see, harmonic measures of certain random walks on cusped graphs of a class of relatively hyperbolic groups.

\subsection{Invariant random walks on graphs}\label{gen_walk}
We refer to Woess \cite{Woess} for the basic theory of random walks on graphs. We mention some of the notions needed in the paper. For a graph $Y$ and a group $\Gamma$ acting on $Y$, a stochastic matrix $P=(p(x,y))_{x,y\in Y}$, that is a matrix such that for all $x\in Y$, $\sum_{y\in Y} p(x,y)=1$, is $\Gamma$-invariant if $$p(gx,gy)=p(x,y),\quad \text{for all} \;x,y\in Y,\; g\in\Gamma.$$ Write $\N_0=\N\cup\{0\}$. A Markov chain with transition probability $P$, is a stochastic process $\{X_i:(Y^{\N_0},\p)\to Y\}_i$, denoted $(Y,P)$, for a probability measure $\p$ on $Y^{\N_0}$, such that for any $n\in\N$ and $x_0,\ldots,x_{n+1}\in Y$ $$\p[X_{n+1}=x_{n+1}|X_n=x_n, \ldots, X_0=x_0]=\p[X_{n+1}=x_{n+1}|X_n=x_n]=p(x_n,x_{n+1}).$$ By the Kolmogorov extension theorem, given a stochastic matrix $P$, there exist corresponding Markov chains $(Y,P)$ with transition probability $P$.

In this paper we work with a $\Gamma$-invariant reversible Markov chain $(Y,P,m)$, which is a $\Gamma$-invariant Markov chain $(Y,P)$, with the property that there is a function, $m:Y\to (0,\infty)$, such that $$m(x)\cdot p(x,y)=m(y)\cdot p(y,x)\quad\text{and}\quad m(gy)=m(y),$$ for all $x,y\in Y$, $g\in\Gamma$.

Given a Markov chain $(Y,P)$, a function $f:Y\to \R$ is $P$-\textit{harmonic} if $$Pf(x):=\sum_y p(x,y)\cdot f(y)=f(x),$$ for all $x\in Y$. A function
$f$ is $P$-\textit{superharmonic} if $Pf\leq f$. The Green function for $(Y,P)$ is defined $G(x,y)=\sum_{n=0}^{\infty} p^{(n)}(x,y)$, for $x,y\in Y$ (here $p^{(0)}(x,y)=\delta_x(y)$). The random walk we consider will be irreducible: for all $x,y\in Y$, there exists $n_{xy}\in\N$, such that $p^{(n_{xy})}(x,y)>0$. 
The walk is recurrent (transient) if there exist $x,y\in Y$, such that $G(x,y)=\infty$ ($G(x,y)<\infty$). An irreducible Markov chain $(Y,P)$ is recurrent if and only if all non-negative superharmonic functions are constants. If the spectral radius of $P$; $\rho(P)=\limsup_{n\to\infty}p^{(n)}(x,x)<1$, then $(Y,P)$ is transient; this will be the case for us.

We define for $x,y\in Y$, $F(x,y)=\sum_{n=0}^\infty\p_x[h_y=n]$, where $h_y=\min\{n\geq 0: X_n=y\}$, as the probability of going from $x$ to $y$, where $\p_x$ is the Markov probability measure for $P$ with initial distribution $\delta_x$. 

Given $A\subset Y$, define $F^A(x,y)=\sum_{n=0}^\infty \p_x[h_A=n, X_n=y]$, where $h_A=\min\{n\geq 0: X_n\in A\}$.
We will need the following inequality got by conditioning on the first visit to $A$ (see \cite[p. 288]{Woess}): $$G(x,y)\geq \sum_{z\in A} F^A(x,z)\cdot G(z,y),$$ for all $x,y\in Y$, with equality if $\p_x$-a.e path from $x$ to $y$ contains an element of $A$ before $y$. 

\subsection{The Martin boundary of a random walk} The reference for this section is again \cite[p. 256]{Woess}.
The Martin kernel of a Markov chain $(Y,P)$ is defined $$K_{y}(x,z)=\frac{G(x,z)}{G(y,z)},\quad\text{for}\;x,y,z\in Y.$$ For $y\in Y$ fixed, the map $z\mapsto K_y(\cdot,z)$ is an embedding from $X$ to $\mathcal{S}_y$, the space of positive superharmonic functions, assuming one at the point $y$ (indeed, $\delta_{z_1}(x)=-\Laplace_x G(x,z_1)=-\Laplace_x G(x,z_2)=\delta_{z_2}(x)$ for all $x\in Y\implies z_1=z_2$, where $\Laplace=P-Id$ is the $P$-Laplacian). The Martin compactification $Y\cup \partial_M Y$ is the closure of $Y\xhookrightarrow{} \mathcal{S}_y$ with respect to the compact-open topology in $\mathcal{S}_y$, and the boundary $\partial_M Y$ is called the Martin boundary. For a sequence $\{z_n\}_n\subset Y$, the limit $\xi=\lim_n z_n\in \partial_M Y$, is often denoted $K(\cdot,\xi)=\lim_n K(\cdot,z_n)$.

A minimal harmonic function normalized at $y\in Y$, is a positive harmonic function $f$ such that $f(y)=1$ and if $f\geq f_1$, where $f_1$ is a positive harmonic function then $f/f_1\equiv$ constant. If $f$ is minimal harmonic, normalized at $y$, then there is $\xi\in\partial_M X$, such that $f\equiv K(\cdot,\xi)$. The set of minimal harmonic functions, denoted $\underline{\partial}_M Y$ is a Borel subset of $\partial_M Y$. There is a correspondence between Borel measures on $\underline{\partial}_M Y$ and $P$-harmonic functions in $Y$ given by the Martin representation formula: given $f\geq 0$ harmonic, there is a unique Borel measure $\mu_y^f$ in $\underline{\partial}_M Y$, such that for all $x\in Y$, $$f(x)=\int_{\underline{\partial}_M Y} K_y(x,\xi)\,d\mu_y^f(\xi).$$

The Martin convergence theorem states that for any $y\in Y$, for the Markov chain $(Y,P)$, $\p_y$-a.e path, $X_n\xrightarrow{n\to\infty}X_\infty$, where $X_\infty$ is a $\underline{\partial}_M Y$-valued random variable with distribution $$\p_x[X_\infty\in E]=\int_{E} K_y(x,\xi)\,d\mu_y^1(\xi),$$ for $x\in Y$, $\mu_y^1$ being the measure corresponding to the constant harmonic function one. Writing $\nu_x(\cdot)=\p_x[X_\infty\in \cdot\,]$, we have for all $x,y,z\in Y$,  $$\nu_y=\mu_y^1\quad\text{and}\quad \frac{d\nu_x}{d\nu_z}(\xi)=\frac{K_y(x,\xi)}{K_y(z,\xi)}=K_z(x,\xi),\quad\text{for}\;\;\nu_z\text{-a.e}\;\;\xi\in\underline{\partial}_M Y,$$ 
We call $\nu=\{\nu_y\}_{y\in Y}$ the \textit{harmonic measure} for the random walk $(Y,P)$.

\subsection{Dimension of measures}
Given a metric space $(A,d)$ and $s\geq 0$, the $s$-Hausdorff measure of a set $E\subset A$ is defined $$\mathcal{H}^s(E)=\lim_{\epsilon\to 0}\; \inf\left\lbrace\sum_i \diam_d(E_i)^s\,\middle|\, E\subset \cup_i E_i, \;\diam_d(E_i)\leq \epsilon\right\rbrace.$$
The Hausdorff dimension of a set $E$, is defined $$\dim_{\mathcal{H}}(E)=\inf\,\{t>0\,\mid\, \mathcal{H}^t(E)=0\}.$$
For a Borel measure $\mu$ in $A$, the (Hausdorff) dimension of $\mu$ is defined $$\dim \mu=\inf\,\{\dim_{\mathcal{H}}(E)\,\mid\, E\subset A, \;\mu(A\setminus E)=0\}.$$
For a finite measure $\mu$, $\dim \mu$ (see Mauldin-Szarek-Urba\'nski \cite[Corollary 8.2]{MSU}),  $$\dim \mu = \text{ess}\sup_{\xi\in A} \liminf_{r\to 0}\,\frac{\log (\mu(B_d(\xi,r)))}{\log r}$$
A finite measure $\mu$ in $A$ is \textit{exact dimensional} if $\lim_{r\to 0}\frac{\log(\mu(B(\xi,r)))}{\log r}$ exists $\mu$-a.e $\xi\in A$ and is a constant.

\subsection{Isoperimetric inequality and the spectral radius of the Markov operator}
Given a reversible Markov chain $(Y,P,m)$, define (see \cite[p. 111]{Woess}) $$\ell^2(Y,m)=\left\{ f:Y\to\R\,\middle|\, \sum_{y\in Y} |f(y)|^2 m(y)<\infty\right\}.$$ There is an isomorphism $$\ell^2(Y,m)\longrightarrow \ell^2(Y), \quad f\mapsto f\cdot \sqrt{m},$$ so that the Markov operator has matrix representation $$P_m(y,x)=(P(\delta_x/\sqrt{m})\,,\,\delta_y/\sqrt{m}\,)_{\ell^2(Y,m)}=\sqrt{m(y)}\cdot p(y,x)/\sqrt{m(x)}$$ as an operator on $\ell^2(Y,m)$. So $P$ is symmetric on $\ell^2(Y,m)$ and its operator norm $\|P\|_{\ell^2(Y,m)}$, is also its spectral radius in $\ell^2(Y,m)$, which is the same as the spectral radius $\rho(P)$ of the transition matrix $P$.

Next we discuss the isoperimetric inequality for $(Y,P,m)$. For a set $E\subset Y$, define the $m$-volume of $E$, $$m(E)=\sum_{y\in E} m(y).$$ Its $(P,m)$-boundary area is defined $$\sigma_{P,m}(\partial E)=\sum_{y\in E}\sum_{\{z|(z,y)\in E(Y), z\neq E\}} m(y)\cdot p (y,z). $$
We say that $(Y,P)$ has the \textit{strong isoperimetric inequality} 
 if there is a constant $\kappa>0$, such that for all finite sets $E$, $$m(E)\leq \kappa\cdot \sigma_{P,m}(\partial E).$$
We use the following result; see \cite[Theorem 10.3]{Woess}.
\bthm \label{10.3}The following are equivalent for $(Y,P,m)$:
\benum
\item $(Y,P,m)$ has the strong isoperimetric inequality.
\item $\rho(P)<1$.
\eenum
\ethm
 Given a $\Gamma$-invariant Markov chain $(\Gamma,P)$ on a discrete countable group, we can define $\mu(g)=p(e,g)$ and consider a random walk on $\Gamma$ with iid increments of law $\mu$. We also use the following \cite{Kesten}.
 
\bthm[Kesten's criterion]\label{Kesten} If \,$\Gamma$ is a non-amenable discrete group, then for a random walk $(\Gamma, \mu)$ driven by i.i.d. increments with law $\mu$, if the support of $\mu$ generates $\Gamma$, then $\mu^{(n)}(e)\leq e^{-\alpha\cdot n}$, for some $\alpha>0$ depending only on $\Gamma, \mu$.
\ethm
So, if $P$ is $\Gamma$-invariant (and irreducible as always), for the random walk $(\Gamma,P)$, we have $\rho(P)<1$.

\bdefn[m-node]\label{mnode} Let $(X,P,m)$ be a reversible random walk. Let $x\in X$. We call the set $N_x:= m^{-1}(m(x))\cap B_X(x,1)$ the $m$-node of $x$. 
\edefn
The following lemma shows that in our situation bounds on $\|P^{2n}\|_{\ell^1\to\ell^\infty}$ (see \cite[Lemma 14.2]{Woess}) follow from $\|P^{2n}\|_{\ell^2\to\ell^2}$ bounds.
\blem\label{2-2-8-1} Let $(X,P,m)$ be a reversible random walk such that $\rho(P)\leq \delta$, and for a constant $C>0$ the following hold:
\benum
\item doubling: for every $x\in X$ and $y\in N_x$, $$\#(N_x\cap N_y)\, \geq \,\frac{C}{m(x)},$$
\vspace{0.2cm}
\item (quasi) invariance: $C\cdot \underset{y\in N_x}{\min}\,p^{(n)}(x,y)\geq \underset{w\in N_v}{\min}\,p^{(n)}(v,w)$, for all $x$, $v$ picked from a level set of $m$ (of the form $m^{-1}(\{c\})$ for $c\in (0,\infty)$) and $n\in\N$,
\vspace{0.2cm}
\item controlled jumps: $p^{(2)}(y,x)\leq C\cdot p^{(6)}(y,z)$, for all $y\in X$ and $x,z\in X$ such that $z\in N_x$.
\eenum
Then, $$\frac{p^{(2n)}(x,x)}{m(x)}\leq  c\cdot \delta^{2n},$$ for a constant $c=c(P,m,C)$, and $ n\in\N$.
\elem 

\bproof
For $v\in X$, for each $n\in\N$, associate a point $v_n\in N_v$ to $v$, which minimizes $p^{(2n)}(v,y)/m(v)$, for $y\in N_v$.

For $y\in N_x$, write $N_x(y)$ for $N_x\cap N_y$.
For $n\in\N$, set $$u_n=\frac{1}{m(x)}\cdot \underset{y\in N_x}{\min}\,p^{(2n)}(x,y).$$ Then for $r,s\in\N$,
$$
\frac{1}{m(x)}\cdot p^{(2(r+s))}(x,x_{r+s})\,\geq\,  \sum_{y\in N_x(x_{r+s})} p^{(2r)}(x,y)\cdot \frac{1}{m(x)}\cdot p^{(2s)}(y,x_{r+s}),
$$ and thus by assumption (2) for some $c_1>0$ depending on $(P, m, C)$, $$u_{r+s}\geq c_1\cdot u_r\cdot u_s.$$ 
Since $\rho(P)<\delta$, $$u_n=\frac{1}{m(x)}\cdot p^{(2n)}(x,x_n)\,\leq\, \frac{1}{m(x)}\cdot\delta^{2n},$$ and thus by Fekete's lemma, $$\lim_{n\to\infty} (u_n)^{1/n}\leq \delta^2.$$ and $$\frac{p^{(2n)}(x,x_n)}{m(x)}\leq c_2\cdot \delta^{2n},$$ for a constant $c_2=c_2(c_1,C)>0$ for all $n\in\N$. 

Next note that by assumption (3), for $n\geq 2$, \beqn\label{firstof2281}\begin{split}p^{(2n)}(x,x) &\leq \sum_y p^{(2n-2)}(x,y)\cdot p^{(2)}(y,x)\\& \leq C\cdot \sum_y p^{(2n-2)}(x,y)\cdot p^{(6)}(y,x_{n+2}) \leq C\cdot p^{(2n+4)}(x,x_{n+2}),\end{split}\eeqn
from which the lemma follows for $n\geq 2$. A similar bound for $n=1$ follows from the third hypothesis.
\eproof

\section{A Random walk on a cusped graph}
\label{construction}
\subsection{A cusped graph (Groves-Manning)}
Let $\Gamma$ be a finitely generated group hyperbolic relative to finitely generated infinite subgroups $H_1,\ldots,H_l$. Let S be a finite, symmetric generating set such that $S\cap H_i$, is a symmetric generating set for $H_i$, $i\in\{1,\ldots, l\}$. We will consider a graph with vertex set $$V_\Gamma=(\Gamma\times\{0\})\,\bigsqcup_{i=1}^l \left(\bigcup_{g\in \Gamma} gH_i\times \N\right).$$
The edge set $E_\Gamma\subset V_\Gamma\times V_\Gamma$ contains pairs of the following types:
\benum
\item $((g,0),(g',0))\in E_\Gamma$ if and only if $(g,g')$ is an edge in $Cay(\Gamma,S)$,
\vspace{0.2cm}
\item for $n\in\N_0$, $g\in\Gamma$ and $h, h'\in H_i$ for some $i\in\{1,\ldots, l\}$, $((gh,n),(gh',n))\in E_\Gamma$ if and only if $d_{(\Gamma,S)}(h,h')\leq a^{n}$.
\vspace{0.2cm}
\item for $n\in\N_0$, $g\in\Gamma$ and $h, h'\in H_i$ for some $i\in\{1,\ldots, l\}$, $((gh,n),(gh',n+1))\in E_\Gamma$ if and only if $d_{(\Gamma,S)}(h,h')\leq a^{n+1}$.
\vspace{0.3cm}
\eenum 
Sometimes to mean $(g,0)\in \Gamma\times\{0\}\subset X_\Gamma$, we only write $g\in \Gamma\subset X$. \vspace{0.2cm}

\underline{Depth function:} The function $y\mapsto n_y$ which assigns the $\Gamma$-orbit of $y\in X$ (the $\N$-coordinate of $x$) will be called the depth function below. The depth is defined to be zero for points in $\Gamma\subset X$.
\vspace{0.2cm}

\underline{Horoballs:} We will call the complete subgraphs of $X$ with vertex subsest $gH_i\times\N$ denoted $(H_i)_g$, \textit{horoballs}. 
\vspace{0.2cm}

\underline{Left-invariant:} Equip the graph $X=X_\Gamma=(V_\Gamma,E_\Gamma)$ with the graph metric, denoted $\rho_X$. It is left-invariant with respect to the following $\Gamma$ action on $X$: $$(g,(g'h,n)):=(gg'h,n),$$ for $g\in\Gamma$, $(g'h,n)\in (H)_g$, $n\in\N_0$. The action extends to an action on the $1$-simplex $X$. Note that $\Gamma$ acts on $X$ properly by isometries, however, the action is not transitive and not acylindrical.
\vspace{0.2cm}

\underline{Hyperbolic:} $(X,\rho_X)$ is a $\delta_X=\delta_X(a,\Gamma)>0$ hyperbolic, proper, geodesic space. This can be seen for example, via its quasiisometry with Groves-Manning's cusped graph. For a detailed description of geodesics in $X_\Gamma$ we refer to \cite[Section 5]{GM}; we state the properties we will use here. The paths described below are called preferred paths in \cite{GM} and this terminology is adopted here.
\vspace{0.2cm}
\benum

\item If $x,y\in (H_i)_g$, consider the points $x_0, y_0\in gH_i$ such that $x=(x_0,n_x)$ and $y=(y_0,n_y)$. Suppose first that $x_0=y_0$. Then denote by $\sigma_{i,g}(x,y)$ the vertical segment with end point $x$ and $y$ in $\{x_0\}\times \N\subset (H_i)_g$. In case $x_0\neq y_0$, denote by $\sigma_{i,g}(x_0,y_0)$ the union of two vertical segments and one horizontal edge, where the depth of the endpoints of the horizontal edge is minimal such, that the points at this depth in the vertical rays\, $\{x_0\}\times\N\subset (H_i)_g$ and $\{y_0\}\times\N\subset (H_i)_g$ can be joined by an edge in $X$. The geodesics joining $x$ and $y$ are then within bounded Hausdorff distance and bounded lengths of the subsegment between $x$ and $y$, denoted $\sigma_{i,g}(x,y)$ of the segment $\sigma_{i,g}(x_0,y_0)$.
\vspace{0.2cm}
\item For $x,y\in X\cup\partial X$ arbitrary, there exists a parametrized (by a connected subset of $\R$) path $p_{xy}$ without backtracks (the same vertex being visited more than once) joining $x$ to $y$ and a collection of ordered horoballs $\mathcal{H}_{xy}=\{A_1,\ldots,A_l\}$ whose union contains $p_{xy}\cap (X\setminus\Gamma)$, such that $p_{xy}\cap A_i$ is of the form $\sigma_{A_i}(x_{i,1},x_{i,2})$ as above, where $x_{i,1}$ is the first point of $p_{xy}$ in $A_i$ and $x_{i,2}$ is the last point of $p_{xy}$ in $A_i$; the path $p_{xy}$ has the property that any geodesic joining $x$ to $y$ is within bounded distance and bounded length of $p_{xy}$.

\eenum

\vspace{0.2cm}
\underline{Visual:} From the combinatorial structure of $X_\Gamma$ and hyperbolicity, we have the following: For $x\in X\cup\partial X$, $y\in X$ there exists $\xi\in\partial X$, such that every geodesic ray $[x,\xi)$, intersects a $\delta_X$-neighbourhood of $y$. Indeed, assuming first that $y\in\Gamma\subset X$; let $(\xi_1,\xi_2)$ be a geodesic line in $X$ joining two points $\xi_1,\xi_2\in\partial X$, and note that it intersects $\Gamma$, and by translation we may assume it passes through $y$. Then by $\delta_X$-thinness of any geodesic triangle $[x,\xi_1,\xi_2]$, either $[x,\xi_1)$ or $[x,\xi_2)$ passes through a $\delta_X$-neighbourhood of $y$. If $y\in (H_i)_g$ for some horoball, then any geodesic line with an end point $\xi_{gH_i}\in\partial X$, the base of $(H_i)_g$, intersects the horosphere containing $y$, and by translation we may assume it passes through $y$. Then the claim follows as in the previous case.
\vspace{0.2cm}

From here on we deal with the case $\Gamma$ has only one virtually nilpotent parabolic subgroup $H$. The general case is only notationally more complex.
\subsection{Construction of a walk}
For any $a>1$, we want $m_0\in (0,\infty)$, $p\in (1,\infty)$, positive sequences $\{p_n\}_{n\in\N_0}$, $\{q_n\}_{n\in\N_0}$ and $\{c_n\}_{n\in\N}$, depending on $a$, such that the following hold:
\benum
\vspace{0.3cm}
\item Normalized: $$\frac{1}{p}+ \frac{g_H(a)}{p_0\cdot g_H(a^{2})}=1,$$
\vspace{0.1cm}
\item and for $n\in\N$,
$$\frac{1}{q_{n-1}}+\frac{g_H(a^{n+1})}{p_n\cdot g_H(a^{n+2})}\leq 1.$$
\vspace{0.1cm}
\item Reversible: $$m_0\cdot \frac{1}{p_0}\cdot \frac{1}{g_H(a^2)}=\frac{c_1}{g_H(a)}\cdot\frac{1}{q_0\cdot g_H(a)},$$
\vspace{0.1cm}
\item and for $n\in\N$, $$\frac{c_n}{g_H(a^n)}\cdot \frac{1}{p_n\cdot g_H(a^{n+2})} = \frac{c_{n+1}}{g_H(a^{n+1})}\cdot \frac{1}{q_n\cdot g_H(a^{n+1})}.$$ 
\vspace{0.1cm}
\item Bounded: There exists $\alpha\geq 1$, such that for all $n\in \N$, $$\alpha\inv \leq c_{n}\leq \alpha, \quad 1 < p_{n}\leq \alpha,\quad 1< q_{n}\leq \alpha.$$
\vspace{0.1cm}
\item There is $0<\beta<1$, such that for all $n\in\N$, $$\frac{q_{n-1}\cdot c_{n+1}\cdot g_H(a^n)}{q_{n}\cdot c_n\cdot g_H(a^{n+1})}<\beta.$$
\vspace{0.3cm}
\eenum
Solutions to above conditions exist:
\benum
\item Set $c_n=1$ for all $n\in\N$.
\item Choose any $p>1$ and $p_0>1$ subject to (1).
\item Choose $\{q_n\}_{n\in\N_0}$, to be a constant sequence, with value larger than $$C\cdot \left(1+\underset{k\in\N}{\max}\,\frac{g_H(a^k)}{g_H(a^{k+1})}\right),$$ for the constant $C>0$ to be fixed next.
\item Choose $p_n$ according to requirement (4), and choose $C$ large enough that $p_n>2$, for all $n\in\N$.
\item Use $q_0$, $p_0$ as chosen above, to choose $m_0$ according to the third requirement (3).
\eenum
\vspace{0.3cm}
Write $d$ for the graph distance in $Cay(\Gamma,S)$. With data as above, define $m:X\to (0,\infty)$ as follows: 
\benum
\vspace{0.2cm}
\item $m((g,0))=m_0$
\vspace{0.2cm}
\item $m((g,n))=\frac{c_n}{g_H(a^n)}\quad\text{for}\;n\in\N,$
\vspace{0.2cm}
\eenum
and $P=(p(x,y))_{x,y\in X}$ as follows:
\benum
\item $\text{for}\; g\sim g'\; \text{in} \;Cay(\Gamma)$, $$p((g,0),(g',0))=\frac{1}{p}\cdot \frac{1}{\text{deg}(Cay(\Gamma))},$$
\item $\text{for}\; g'\in gH\;\text{where}\; \di(g,g')\leq a^{n+1}\; \text{in} \;Cay(\Gamma)\; \text{and}\;n\geq 0,$
$$p((g,n+1),(g',n))=\frac{1}{q_n\cdot g_H(a^{n+1})}.$$
\item $\text{for}\; g'\in gH\;\text{where}\; \di(g,g')\leq a^{n+1}\; \text{in} \;Cay(\Gamma)\; \text{and}\;n\geq 0,$ $$p((g,n),(g',n+1))=\frac{1}{p_n\cdot g_H(a^{n+2})}.$$ 

\item for $g'\in gH$, $n\in\N$ the remaining probability for jump starting at $(g,n)$ is distributed among neighbours of $(g,n)$ in $m\inv(\{m((g,n))\})$ in any $\Gamma$-invariant way (for example, uniformly). 

\eenum
Then $m$ is an invariant measure for $P$ and $(X,P,m)$ defines a $\Gamma$-invariant reversible Markov chain. 
\subsection{The $\rho(P)$ bound} We next check the strong isoperimetric inequality for $(X, P, m)$.
\blem[Strong isoperimetric inequality]\label{strongIS} Let $A$ be finite subset of $X$. There exists $\kappa>0$, such that $$m(A)\,\leq \,\kappa\cdot  \sigma_{P,m}(\partial A).$$
\elem
\bproof
We write $\sigma=\sigma_{P,m}$. Let $H_\alpha=g_\alpha H$, $\alpha\in\mathcal{J}\subset\N$ be the cosets of $H$. We may assume that $A$ as a subspace of the $1$-simplex $X$ is connected since $P$ is nearest-neighbour.  
Set $$A_i^\alpha=A\cap (H_\alpha\times\{i\}),$$ for $\alpha\in\mathcal{J}$ and $n\in\N$ and $A_0=A\cap(\Gamma\times\{0\})$. 

\textbf{Step 1} \underline{Compact part}:
A computation shows $$\sigma(\{(x,y)\mid x\in A_0, y\in (\Gamma\times\{0\})\setminus A_0\}) \,\geq\, C_1\cdot \sigma_{Cay(\Gamma)}(\partial_{Cay(\Gamma)} A_0),$$ where $C_1=C_1(p,p_0,f,a)>0$ is a constant and $\sigma_{Cay(\Gamma)}$ is the boundary area for the simple random walk in $\Gamma$. Then, by Kesten's theorem (Theorem \ref{Kesten}) and Theorem \ref{10.3}, there is $\kappa_0>0$, such that
\beqn\label{firstofIS}
m_0\cdot \#A_0\, \leq \,\kappa_0\cdot\sigma(\{(x,y)\mid x\in A_0, y\in (\Gamma\times\{0\})\setminus A_0\}).
\eeqn

\textbf{Step 2} \underline{Cusp part}:
Next, let $n_\alpha\in\N$ be the smallest $i\in\N$ such that $A_i^\alpha=\emptyset$. For $n\in\N_0$, and $x\in A^\alpha_{n+1}$ set $$B_n(x)=\{y\in A^\alpha_n\,\mid\,x\sim y\},$$ and for $n\in\N$, $x\in A^\alpha_{n}$, set $$\tilde{B}_n(x)=\{y\in A^\alpha_{n+1}\,\mid\,x\sim y\}.$$
By connectedness of $A$, $A_{i+1}^\alpha=\emptyset$ for all $i\,>\, n_\alpha$. 
Then we have
\beqn\label{secondofIS}
\begin{split}
\sigma(\partial\, (\bigcup_{i=1}^{n_\alpha}A_i^\alpha))&= \frac{c_1}{g_H(a)}\cdot\frac{1}{q_0\cdot g_H(a)}\sum_{x\in A_1^\alpha}g_H(a) \,+\, \frac{c_1}{g_H(a)}\cdot\frac{1}{p_1\cdot g_H(a^3)}\cdot \sum_{x\in A_1^\alpha}(g_H(a^3)-\#\tilde{B}_2(x))\\& + \frac{c_2}{g_H(a^2)}\cdot \frac{1}{q_1\cdot g_H(a^2)}\sum_{x\in A_2^\alpha}(g_H(a^2)-\#B_1(x)) \,+\,\dots \\& + 
\frac{c_{n_\alpha-1}}{g_H(a^{n_\alpha-1})}\cdot \frac{1}{p_{n_\alpha-1}\cdot g_H(a^{n_\alpha+1})}\sum_{x\in A_{n_\alpha-1}}(g_H(a^{n_\alpha+1})-\#\tilde{B}_{n_\alpha}(x)) \\& + 
\frac{c_{n_\alpha}}{g_H(a^{n_\alpha})}\cdot \frac{1}{q_{n_\alpha-1}\cdot g_H(a^{n_\alpha})} \sum_{x\in A_{n_\alpha}^\alpha}(g_H(a^{n_\alpha})-\#B_{n_\alpha-1}(x)) \\& + \frac{c_{n_\alpha}}{g_H(a^{n_\alpha})}\cdot \frac{1}{p_{n_\alpha}\cdot g_H(a^{n_\alpha+2})}\cdot \sum_{x\in A_{n_\alpha}^\alpha}g_H(a^{n_\alpha+1}) \\&= 
\frac{c_1}{g_H(a)}\cdot\frac{1}{q_0}\cdot \#A_1^\alpha \,+\, \sum_{i=1}^{n_\alpha-1}\frac{c_i}{g_H(a^i)}\cdot\frac{1}{p_i\cdot g_H(a^{i+2})}\sum_{x\in A_i^\alpha}(g_H(a^{i+1})-\#\tilde{B}_{i+1}(x)) \\& + 
\sum_{i=1}^{n_\alpha-1}\frac{c_{i+1}}{g_H(a^{i+1})}\cdot\frac{1}{q_{i}\cdot g_H(a^{i+1})}\sum_{x\in A_{i+1}^\alpha}(g_H(a^{i+1})-\#B_{i}(x)) \\& +\frac{c_{n_\alpha}}{g_H(a^{n_\alpha})}\cdot \frac{g_H(a^{n_\alpha+1})}{p_{n_\alpha}\cdot g_H(a^{n_\alpha+2})}\cdot \#A_{n_\alpha}^\alpha.
\end{split}
\eeqn

Note that for $i\in \{1,\ldots, n_i-1\}$,
\beqnn\begin{split}\frac{1}{g_H(a^{i+1})}\sum_{x\in A_i^\alpha}\#\tilde{B}_{i+1}(x)&=\frac{1}{g_H(a^{i+1})}\sum_{x\in A_i^\alpha}\sum_{y\in A_{i+1}^\alpha}\ind_{E_\Gamma}(x,y) \\& = \sum_{y\in A_{i+1}^\alpha}\frac{\sum_{x\in A_i^\alpha}\ind_{E_\Gamma}(x,y)}{g_H(a^{i+1})} \\& \leq \#A_{i+1}^\alpha, \end{split}\eeqnn
and similarly,
$$\frac{1}{g_H(a^{i+1})}\sum_{x\in A_{i+1}^\alpha}\#B_{i}(x)\leq \#A_i^\alpha.$$
It follows that 
$$\frac{1}{g_H(a^{i+1})}\sum_{x\in A_i^\alpha}(g_H(a^{i+1})-\#\tilde{B}_{i+1}(x))\geq \max\{0, \#A_i^\alpha-\#A_{i+1}^\alpha\},$$ and
$$\frac{1}{g_H(a^{i+1})}\sum_{x\in A_{i+1}^\alpha}(g_H(a^{i+1})-\#B_{i}(x))\geq \max\{0, \#A_{i+1}^\alpha-\#A_{i}^\alpha\},$$

Then from \eqref{secondofIS} and property (4) of $P$ (reversibility), and property (6),
\beqn\label{thirdofIS}
\begin{split}
\sigma(\partial\, (\bigcup_{i=1}^{n_\alpha}A_i^\alpha)) &\geq \frac{c_1}{g_H(a)}\cdot\frac{1}{q_0}\cdot\#A_1^\alpha\,+ \, \sum_{i=1}^{n_\alpha-1}\frac{c_{i+1}}{g_H(a^{i+1})}\cdot\frac{1}{q_i}\,|\#A_{i+1}^\alpha-\#A_i^\alpha| + \\& +\frac{c_{n_\alpha+1}}{g_H(a^{n_\alpha+1})}\cdot \frac{1}{q_{n_\alpha}}\cdot \#A_{n_\alpha}^\alpha \\&\geq 
\sum_{i=1}^{n_\alpha} \left(\frac{c_i}{q^{i-1}\cdot g_H(a^{i})}-\frac{c_{i+1}}{q_i\cdot g_H(a^{i+1})}\right) \cdot \#A_i^\alpha \\& \geq  \sum_{i=1}^{n_\alpha} \frac{1}{q_{i-1}}\cdot \left(1-\frac{q_{i-1}\cdot c_{i+1}\cdot g_H(a^i)}{q_{i}\cdot c_i\cdot g_H(a^{i+1})}\right)\cdot \frac{c_i}{g_H(a^i)}\cdot \#A_i^\alpha \\&\geq \frac{1}{\overline{\kappa}}\cdot m(\bigcup_{i=1}^{n_\alpha}A_i^\alpha),
\end{split}
\eeqn
where $\overline{\kappa}$ depends only on $g_H$, $a$, and the data associated to $P$ and $m$ (and not on $A$). 

\textbf{Step 3} \underline{Summing}: 
Next, note that by definition of $\sigma$,
\beqn\label{fourthofIS}
\begin{split}
\sigma(\partial A) &\geq  \sigma(\{(x,y)\mid x\in A_0, y\in (\Gamma\times\{0\})\setminus A_0\}) \\& + 
\sum_\alpha \left(\sigma(\partial\, (\bigcup_{i=1}^{n_\alpha}A_i^\alpha))\,-\, \frac{c_1}{g_H(a)}\cdot\frac{1}{q_0\cdot g_H(a)}\sum_{x\in A_1^\alpha}\#B_0(x) \right).
\end{split}
\eeqn
Write $$\mathcal{I}:=\left\{\alpha: \frac{1}{\overline{\kappa}}\cdot m(\bigcup_{i=1}^{n_\alpha}A_i^\alpha)\leq 2\cdot \frac{c_1}{g_H(a)}\cdot\frac{1}{q_0\cdot g_H(a)}\sum_{x\in A_1^\alpha}\#B_0(x)\right\}.$$
Then by \eqref{firstofIS}, \eqref{thirdofIS}, and \eqref{fourthofIS},
\beqnn
\begin{split}
m(A) & = m_0\cdot \#A_0 \,+\, \sum_{\alpha\in\mathcal{I}} m(\bigcup_{i=1}^{n_\alpha}A_i^{\alpha})\,+\, \sum_{\alpha\notin\mathcal{I}} m(\bigcup_{i=1}^{n_\alpha}A_i^{\alpha}) \\& \leq 
\,\kappa_0\cdot\sigma(\{(x,y)\mid x\in A_0, y\in (\Gamma\times\{0\})\setminus A_0\}) \\& + 2\cdot\overline{\kappa}\cdot \frac{c_1}{g_H(a)}\cdot\frac{1}{q_0\cdot g_H(a)}\sum_{\alpha\in\mathcal{I}}\sum_{x\in A_1^\alpha}\#B_0(x) \\& + 
2\cdot\overline{\kappa}\cdot\sum_{\alpha\notin\mathcal{I}} \left(\frac{1}{\overline{\kappa}}\cdot m(\bigcup_{i=1}^{n_\alpha}A_i^\alpha)\,-\,\frac{c_1}{g_H(a)}\cdot\frac{1}{q_0\cdot g_H(a)}\sum_{x\in A_1^\alpha}\#B_0(x)\right) \\& \leq \,
\kappa_0\cdot\sigma(\{(x,y)\mid x\in A_0, y\in (\Gamma\times\{0\})\setminus A_0\}) \, +\,
\frac{2\cdot \overline{\kappa}\cdot c_1}{q_0\cdot g_H(a)}\#A_0 \\& + 
2\cdot\overline{\kappa}\cdot\sum_{\alpha\notin\mathcal{I}}\left(\sigma(\partial\, (\bigcup_{i=1}^{n_\alpha}A_i^\alpha))\,-\, \frac{c_1}{g_H(a)}\cdot\frac{1}{q_0\cdot g_H(a)}\sum_{x\in A_1^\alpha}\#B_0(x) \right) \\& \leq \kappa_0\cdot\left(1+\frac{2\cdot\overline{\kappa}\cdot c_1}{q_0\cdot g_H(a)\cdot m_0}\right)\cdot \sigma(\{(x,y)\mid x\in A_0, y\in (\Gamma\times\{0\})\setminus A_0\}) \\& + 2\cdot\overline{\kappa}\cdot\sum_{\alpha}\left(\sigma(\partial\, (\bigcup_{i=1}^{n_\alpha}A_i^\alpha))\,-\, \frac{c_1}{g_H(a)}\cdot\frac{1}{q_0\cdot g_H(a)}\sum_{x\in A_1^\alpha}\#B_0(x) \right) \\& \leq 
2\cdot \max\left\{2\cdot\overline{\kappa},\,\kappa_0\cdot\left(1+\frac{2\cdot\overline{\kappa}\cdot g_H(a)}{p\cdot deg(Cay(\Gamma))}\right)\right\}\cdot \sigma(\partial A).
\end{split}
\eeqnn
This concludes the proof.
\eproof

\section{The Bowditch boundary as a Martin boundary}\label{BowditchisMartin}

Consider the reversible Markov chain $(X,P,m)$ introduced in \S \ref{construction}. In this section we identify its Martin boundary with the Bowditch boundary of $(\Gamma,H)$. We start by noting some properties of the walk.
\blem\label{replacement}
Let $x, y, z\in X$ and assume that $z\in N_x$, the $m$-node of $x$. Then for the reversible Markov chain $(X,P,m)$, $$p^{(2)}(y,x)\leq C_1\cdot p^{(6)}(y,z),$$ for a constant $C_1\geq 1$.
\elem
\bproof
We may assume $y\in B_X(x,2)$. Since $m(x)=m(z)$, by reversibility of $P$, it suffices to check that $p^{(2)}(x,y)\leq C_1\cdot p^{(6)}(z,y)$, for some $C_1\geq 1$. This is clear in case $x\in \Gamma\times\{0\}$. Then for $x$ in a horoball one notices that $\underset{y'\in X}{\max} \;p^{(2)}(x,y')$ is within constant multiples of $m(x)$. Let us assume first $n_y=n_x$. Then paths $\{x_0=z,x_1,x_2,x_3,x_4,x_5=y\}$, where the even-labelled vertices are in the horosphere containing $x$, and the odd-labelled vertices are in the horosphere $n_x+1$, can be chosen from $z$ to $y$. Each of the vertices $x_i$ for $i\in\{1,2,3,4\}$ can be chosen in at least $c/m(x)$ ways, for some $0<c<1$ such that the transition probability from $x_i$ to $x_{i+1}$, for $i\in\{0,\ldots,4\}$ is at least $c'\cdot m(x)$ for some $c'>0$ (depending on data associated to $P$ and $m$ and using polynomial growth $g_H$ of $H$). The cases $n_y\in\{n_x-2,n_x-1,n_x+1,n_x+2\}$ are similar.
\eproof

Next is an `uniform irreducibility' property.
\blem\label{basicQI}
There exists $\lambda\in(0,1)$ and $C>0$, such that for $x,y\in X$, $$\frac{p^{(r(x,y))}(x,y)}{m(y)}\geq C\cdot \lambda^{\rho_X(x,y)},$$ where $0\leq r(x,y)-\rho_X(x,y)\leq c_0,$ for a constant $c_0>0$. \elem
\bproof
First consider the case that $x,y$ both lie in the same horoball. Let $n_x$ and $n_y$ be the depths of the horocycles containing $x$ and $y$ respectively. We may assume without loss of generality that $n_x\geq n_y\geq 0$. Let $k_{xy}\in \N$ be such that there is a preferred path $p_{xy}$ which passes through a point $z=z_{xy}$ in the $k_{xy}$-th horosphere. Write $r_x=\rho_X(x,z)$ and $r_y=\rho_X(y,z)$. For $z' \in N_z$ a path $\{x_1=x,\ldots x_{r_x}=z'\}$ of length $r_x$ from $x$ to $z'$, there are at least $c/m(x_i)$ ways to choose $x_i$ such that the transition probability from $x_i$ to $x_{i+1}$ is $c'\cdot m(x_{i+1})$, for some $c,c'>0$ for $i\in\{2,\ldots,r_x-1\}$ (relying again on the polynomial growth of $H$). Then $$p^{(r_x)}(x,z')\geq \frac{C'\cdot m(n_x+r_x)}{C^{(r_x-1)}},$$ for some constant $C>1$. After a similar consideration for $p^{(r_y)}(y,z'')$, and $z''\in N_z$ we note
\beqn \begin{split}
p^{(r_x+r_y+4)}(x,y) &\geq \sum_{z'\in N_z}p^{(r_x)}(x,z')\cdot p^{(4)}(z',z'')\cdot p^{(r_y)}(z'',y) \\&
\geq \frac{C'}{C^{(r_x-1)}}\cdot c'\cdot m(z'')\cdot p^{(r_y)}(z'',y) \\& \geq  \frac{C'}{C^{(r_x-1)}}\cdot c'\cdot m(y) \cdot p^{(r_y)}(y,z'') \geq \frac{C''}{C^{(r_x+r_y-2)}}\cdot m(y),
\end{split}
\eeqn
for some $C''>0$.

Next consider the case when $x,y\in\Gamma\times\{0\}$ and a preferred path joining them is contained in $\Gamma\times\{0\}$. The claim in this case is clear. In the case when a preferred path joining points $x$ and $y$ intersects horoballs as well as $\Gamma\times\{0\}$, we apply the above discussion to each of the components that lie in the horoballs and in $\Gamma\times\{0\}$ separately. A function $r$ can be defined taking the product of suitable lower bounds for probabilities corresponding to each of the pieces to satisfy the claim.
\eproof

\subsection{A Harnack inequality for $(X,P,m)$}

We start by noting another basic property of the transition matrix $P$.
\blem\label{basic_compar} For $x,y\in X$, $y'\in N_y$, $x\neq y$ and $2\leq n\in \N$, $$p^{(n)}(x,y)\lesssim p^{(n+4)}(x,y').$$ In particular, there is $C>0$ so that, 
$$G(x,y)\leq C\cdot G(x,y'),$$ for any $y'\in N_y$.
\elem
\bproof From Lemma \ref{replacement}, we have for $n\geq 2$,
\beqnn
\begin{split} p^{(n)}(x,y) &= \sum_{v} p^{(n-2)}(x,v)\cdot p^{(2)}(v,y)\\& 
\leq C_1 \cdot \sum_{v} p^{(n-2)}(x,v)\cdot p^{(6)}(v,y')\\&
\leq C_1\cdot p^{(n+4)}(x,y').
\end{split}
\eeqnn
For the second claim, observe that $p(x,y)\neq 0$, then $x\in B_X(y',2)$ and again (similarly to Lemma \ref{replacement}), $\min\{p^{(5)}(x,y'), p^{(6)}(x,y')\}\geq C_2\cdot m(y)\geq C_3\cdot p(x,y)$, for some
$C_2,C_3>0$. Summing we get the claim.
\eproof

By Theorem \ref{10.3} and Lemma \ref{strongIS} we have $p^{(2n)}(x,x)\leq \delta^{2n}$ for some $0<\delta<1$. We need the following stronger property. 
\bcor\label{8to1}
For the reversible Markov chain $(X,P,m)$, there exist numbers $C>0$ and $0<\delta<1$, such that for any $x\in X$ and $n\in \N$, $$\frac{p^{(2n)}(x,x)}{m(x)}\leq \,C\cdot \delta^{2n},$$ and consequently, $$\frac{p^{(n)}(x,y)}{m(y)}\leq C'\cdot \delta^n,$$ for all $x,y\in X$, $n\in\N$ and some $C'>0$.
\ecor 
\bproof
Note that there is $0<\delta<1$ such that $\rho(P)<\delta$ by Lemma \ref{strongIS}. The claim is clear for $x\in\Gamma\times\{0\}$. Suppose then that $x$ lies in a horoball. We observe that the points in $N_x$ correspond to points in a ball of radius $a^{n_x}$ in $H$ as a subgraph of $Cay(\Gamma)$. The first assumption of Lemma \ref{2-2-8-1} is then verified, using polynomial growth of $H$. The second assumption uses invariance of $P$ and $m$ by the $\Gamma$ action and the third assumption is checked in Lemma \ref{replacement}. Then the first part of the corollary follows from Lemma \ref{2-2-8-1}. For the second see \cite[Lemma 14.2]{Woess} and write $P^{2n+1}=P^{2n}\cdot P$.
\eproof

\blem[Harnack inequality for positive superharmonic functions]\label{harnack} The following hold:
\benum
\item\label{harnack1} If $x,y,z\in X$, and $x\neq y$, then $G(x,y)\leq C^{\,\rho_X(x,z)}\cdot G(z,y).$
\vspace{0.2cm}
\item\label{harnack2} For $x,y\in X$, $m(x)\cdot G(x,x)\leq C^{\,\rho_X(x,y)}\cdot G(y,x)$
\vspace{0.2cm}
\item For $u$ positive, $P$-superharmonic, and $x,y\in X$, $m(x)\cdot u(x) \leq C^{\,\rho_X(x,y)}\cdot u(y).$
\vspace{0.2cm}
\eenum
\elem
\bproof
For $n\geq 2$ (recalling $r(\cdot,\cdot)$ from Lemma \ref{basicQI})
\beqn\label{firstofharnack}
\begin{split}
\frac{p^{(n)}(y,x)}{m(x)}\cdot \frac{p^{(r(x,z))}(x,z)}{m(z)}&\leq C_1\cdot \sum_{x'\in N_x}p^{(n+4)}(y,x')\cdot \frac{p^{(r(x,z))}(x',z)}{m(z)} \\& \leq C_1\cdot \frac{p^{(n+r(x,z)+4)}(y,z)}{m(z)},
\end{split}
\eeqn
from Lemma \ref{basic_compar}. Then from Lemma \ref{basicQI},
$$\frac{p^{(n)}(y,x)}{m(x)}\leq \frac{C_1}{C\cdot \lambda^{\rho_X(x,z)}}\cdot \frac{p^{(n+r(x,y)+4)}(y,z)}{m(z)}.$$ A similar bound holds for $n=1$, noting $p(y,x)\leq C'\cdot p^{(6)}(y,x')$, for some $C>0$, as in the proof of Lemma \ref{basic_compar}. Then summing and using reversibility of $P$ by $m$ gives (1). 

The second claim follows from Corollary \ref{8to1}, that is \beqn\label{secondofharnack} 1\leq G(x,x) = p^{(0)}(x,x) + \sum_{n=1}^{\infty} p^{(n)}(x,x)\leq 1+ \sum_{n=1}^{\infty}C\cdot \delta^{n}\cdot m(x)= 1+\frac{C\cdot \delta\cdot m(x)}{1-\delta},\eeqn and Lemma \ref{basicQI}.

For (3), first note that for a function $f\geq 0$, by (1) and (2) 
$$Gf(x)=\sum_y G(x,y)\cdot f(y)\leq \frac{C^{\rho_X(x,z)}}{m(x)}\cdot \sum_y G(z,y)\cdot f(y)=\frac{C^{\rho_X(x,z)}}{m(x)}\cdot Gf(z).$$
The claim for general $u$ follows by approximation by writing $u=Gf+h$ for some function $f\geq 0$ and harmonic function $h\geq 0$ (Riesz decomposition, see \cite[Theorem 24.2]{Woess}) and $h=\underset{n\to\infty}{\lim} Gf_n$ as an increasing limit for some sequence $\{f_n\}_n$ of non-negative functions (\cite[Theorem 24.6]{Woess}). Then $u=\underset{n\to\infty}{\lim}G\phi_n$ as an increasing limit, where $\phi_n$ are non-negative functions.
\eproof
The  Harnack inequality for harmonic functions holds. We will not be using this fact however.
\bcor[Harnack inequality for harmonic function] 
If $u$ is a non-negative $P$-harmonic function, then there exists $C\geq 1$ such that $u(x)\leq C^{\rho_X(x,y)} \cdot u(y)$.
\ecor
\bproof
This follows from \eqref{harnack1} in Lemma \ref{harnack} and the Martin representation formula for harmonic functions.
\eproof
For any set $B\subset X$, define $\partial B=\{z\in X\setminus B\,\mid\, (z,y)\in E_\Gamma(X),\;\text{for some}\; y\in B\}.$

\blem[Triangle inequality] \label{triangle}There exists $C\geq 1$ such that for $x,y,z\in X$, $y\notin\{x,z\}$,
\beqn\label{triangleeq1}\frac{F(x,y)}{m(y)}\cdot F(y,z)\leq C\cdot F(x,z).\eeqn
\elem
\bproof
We have by Lemma \ref{basic_compar}, and $y\neq z$,
\beqn\label{firstoftriangle}
\begin{split}
G(x,z) &\geq \sum_{w\in \partial N_y} F^{\partial N_y}(x,w)\cdot G(w,z)\geq \frac{1}{C} \cdot G(y,z)\cdot \sum_{w\in \partial N_y}F^{\partial N_y}(x,w).
\end{split}
\eeqn
Next if $x\notin \overline{N}_y$,
\beqn\label{secondoftriangle}
\begin{split}
G(x,y)&=\sum_{w\in\partial N_y} F^{\partial N_y}(x,w)\cdot G(w,y)  \leq C'\cdot m(y)\cdot \sum_{w\in\partial N_y}F^{\partial N_y}(x,w)
\end{split} 
\eeqn
(where $C'>0$ comes from Corollary \ref{8to1}). Then the lemma follows comparing \eqref{firstoftriangle} and \eqref{secondoftriangle} in case $x\notin \overline{N}_y$. If however $x\in \overline{N}_y$, pick $x'\in X$ such that $\rho_X(x,x')\leq 2$ and $x'\notin \overline{N}_y$. Then the lemma follows by equations \eqref{firstoftriangle} and \eqref{secondoftriangle} applied to $x',y,z$ and the Harnack inequality \eqref{harnack1} of Lemma \ref{harnack} (where $x\neq y$ is used).
\eproof

\subsection{The Martin boundary of $(X,P)$}
Next, the properties of the reversible Markov chain $(X,P,m)$ gathered so far are employed to identify the Martin boundary of the random walk $(X,P)$ (the proofs of the results stated below, follow by suitably modifying the proofs in \cite[p. 288]{Woess}).

\bthm[Ancona inequality]\label{ancona} Let $(X,P,m)$ be a reversible Markov chain. If:
\benum
\item the graph $X$ is Gromov hyperbolic,
\item Harnack inequalities \eqref{harnack1} and \eqref{harnack2} hold,
\item triangle inequality \eqref{triangleeq1} holds,
\eenum
then there exist constants $C_0,C\geq 1$ such that for $x,y,z\in X$ with $\dist_X(y,[x,z])\leq r$, we have
\beqn\label{anconain}
F(x,z)\leq C_0\cdot C^{2r}\cdot \frac{F(x,y)}{m(y)}\cdot F(y,z),
\eeqn
for any geodesic $[x,z]$ joining $x,z$.
\ethm
\brem
Note from Lemma \ref{triangle} that the factor of inverse $m(\cdot)$ on the right hand side of \eqref{anconain} as stated above is necessary.
\erem

\bthm[Martin boundary]\label{Martin} Let $(X,P,m)$ be a reversible Markov chain. If:
\benum
\item the graph $X$ is Gromov hyperbolic,
\item Harnack inequalities \eqref{harnack1} and \eqref{harnack2} hold,
\item the triangle inequality \eqref{triangleeq1} holds,
\item the Ancona inequality \eqref{anconain} holds,
\eenum
then the minimal Martin boundary and the Martin boundary of the random walk $(X,P)$ are homeomorphic to the visual boundary $\partial X$ of $X$; for any $x\in X$, the embedding $y\mapsto K_x(\cdot,y)$, $y\in X$, extends homeomorphically to $X\cup\partial X\longrightarrow X\cup\partial_M X$.
\ethm
\brem[$r$-Martin boundary]
The steps leading to Theorem \ref{Martin} also apply to the $r$-Martin boundary $\partial_{M,r}X$, for $1\leq r <1/\rho(P)$, with the corresponding $r$-Green function $G_r(x,y)=\sum_{n\in\N_0} r^n\cdot p^{(n)}(x,y)$. Theorem \ref{Martin} then also holds when $1\leq r<1/\rho(P)$.
\erem

\section{The harmonic measure as a conformal density}\label{metricmeasurestructure}
We start by introducing a Green metric on $X$, corresponding to the Markov chain $(X,P,m)$. 
\subsection{A Green metric}\label{green_metric} By Lemmas \ref{8to1}, \ref{triangle} and equation \eqref{secondofharnack}, there exists $C'\geq 1$, such that $$\frac{G(x,z)}{C'\cdot m(z)}\geq \frac{G(x,y)}{C'\cdot m(y)}\cdot \frac{G(y,z)}{C'\cdot m(z)},$$ and $$\frac{G(x,y)}{C'\cdot m(y)} < e^{-1},$$ for all $x,y,z\in X$, where $y\notin\{x,z\}$. Let $C$ be the infimum of all such $C'$ and set for $x,y\in X$, \[G_m(x,y)=\begin{cases}\frac{G(x,y)}{C\cdot m(y)} &x\neq y\\ 1& x=y.\end{cases}\]

\blem For $x,y\in X$, the function
$$\rho_G(x,y)=-\log(G_m(x,y)),$$ is a $\Gamma$-invariant metric.\elem 
The precise choice of the constant $C$ is insignificant for us as a different (admissible) choice leads to a function which is a bounded distance from the function defined above.
\bproof The triangle inequality follows by noting $$G_m(x,y)\cdot G_m(y,z)\leq G_m(x,z),$$ for all $x,y,z\in X$. Symmetry follows by reversibility of $(X,P,m)$. It is clear $$\rho_G(x,y)=0 \iff x=y.$$ Invariance with respect to $\Gamma$ action comes from the $\Gamma$-invariance of $P$ and $m$. \eproof

\bthm[Hyperbolicity and QI]\label{QI} The metrics $\rho_G$ and $\rho_X$ are quasiisometric and $(X,\rho_{G})$ is a hyperbolic quasiruled space.
\ethm 
\bproof
Note that $\rho_X(x,y)=1\implies 1\leq \rho_G(x,y)\leq C_0,$ where $C_0\geq 1$ comes from the second Harnack inequality \eqref{harnack2}. If $x=x_0,\ldots,x_n=y$ is a geodesic path joining $x$ and $y$ in $(X,\rho_X)$, then $$\rho_G(x,y)\leq \sum_{i=0}^{n-1} \rho_G(x_i,x_{i+1})\leq C_0\cdot n=C_0\cdot \rho_X(x,y).$$ By the Harnack inequality \eqref{harnack2} again, $$\rho_X(x,y)\leq C_1\cdot \rho_G(x,y)+C_2,$$ for some $C_1\geq 1, C_2\geq 0$.
Then the geodesics of $(X,\rho_X)$ are quasigeodesics of $(X,\rho_G)$. Ancona inequality (Lemma \ref{anconain}) shows that the geodesics of $(X,\rho_X)$ are quasirulers in $(X,\rho_G)$. Then $(X,\rho_G)$ is quasiruled. It follows from Theorem A.1, \cite{BHM} that $(X,\rho_G)$ is also hyperbolic.
\eproof

\subsection{The Green-Busemann function and conformality of the harmonic measure}\label{QC} For $x,y,z\in X$, set  $$\beta_{z}^{G}(x,y)=\rho_G(x,z)-\rho_G(y,z).$$ For $\{z_n\}_n\subset X$, $z_n\to\xi\in\partial X$, by Theorem \ref{Martin}, $$\beta_{z_n}^G(x,y)=-\log\left(\frac{G(x,z_n)}{G(y,z_n)}\right)=-\log\left(\frac{K(x,z_n)}{K(y,z_n)}\right)\xrightarrow{n\to\xi} -\log\left(\frac{K(x,\xi)}{K(y,\xi)}\right),$$ where $K(\cdot,\cdot)=K_e(\cdot,\cdot)$. 
Then extend to the boundary $\beta^G:\partial X\times X\times X\to\R$, $$\beta^G_{\xi}(x,y)=-\log\left(\frac{K(x,\xi)}{K(y,\xi)}\right),$$
which is the Busemann function of the metric $\rho_G$. For the harmonic measures $\{\nu_x\}_{x\in X}$, we get for $\nu_x$-a.e $\xi\in\partial X$,
\beqnn
\frac{ d\nu_y}{d \nu_x}(\xi)=e^{-\beta^G_\xi(y,x)},
\eeqnn
and consequently, by equivariance of the $\Gamma$ action (on $(X,P,m)$)
\beqn\label{conformality}
\frac{d\gamma_\ast \nu_x}{d \nu_x}(\xi)=e^{-\beta^G_\xi(\gamma x,x)},
\eeqn
for $\nu_x$ a.e $\xi\in\partial X$ and $\gamma\in \Gamma$ (where $\gamma_\ast$ is the push-forward: $\gamma_\ast\nu_x(E)=\nu_x(\gamma^{-1}E)$). Thus the harmonic measure is a $\Gamma$-conformal density of dimension one with respect to the Green metric.

\section{The shadow lemmas}\label{doubling}
We equip the Bowditch boundary with the visual metrics $\{d_{\epsilon,X}\}_{0<\epsilon<\epsilon_1}$, $\{d_{\epsilon,G}\}_{0<\epsilon<\epsilon_2}$, induced respectively, by the graph metric $\rho_X$ and the Green metric $\rho_G$. In this section we discuss the metric properties of the harmonic measure $\nu$. 
\brem[Ergodicity of $(X,\nu,\Gamma)$] \label{ergodicity}Let us note that $\nu$ is non-atomic follows easily: Suppose to the contrary, $\xi\in\partial X$ is an atom for $\nu_e$. Then $$K(x,\xi)\cdot \nu_e(\{\xi\})\leq 1,$$ so we have $K(\cdot,\xi)$ is bounded. Then, $\beta^G_\xi(\cdot,e)$ is bounded below. This is a contradiction (with the hyperbolicity of $(X,\rho_G)$). Thus the set of conical limit points has full $\nu$ measure. From this, the ergodicity of the action of $\Gamma$ with respect to the harmonic measure follows by another classical argument (see \cite[Theorem 4.4.4]{Nic}) given the doubling property of $\nu$ below. If $\mu=\{\mu_x\}_{x\in X}$ is a $D_\Gamma$-dimensional Patterson-Sullivan density for $\Gamma$, then $\mu$ is also doubling and non-atomic (see Remark \ref{imp_1} below), and the same holds.
\erem

\bthm\label{double}
The harmonic measure is Ahlfors-regular for $(\partial X,d_{G,\epsilon_G})$.
\ethm

The proof of Theorem \ref{double} follows from Lemma \ref{shadow} below. 

\blem[Sullivan's shadow lemma]\label{SS} Let $(Y,\rho)$ by a hyperbolic quasiruled space. Suppose a group $G$ acts non-elementarily on $Y$ by isometries. Let $\{\nu_y\}_{y\in \Gamma}$ be a $G$-conformal density of dimension $\beta>0$ in $\partial Y$. Then, for each $y\in Y$ there exists $R_S=R_S(y)>0$, such that for any $r\geq R_S$, there exists $C=C(r,y)\geq 1$ so that for all $g\in G$, $$\frac{1}{C}\cdot e^{-\beta\cdot \rho(y,g y)}\leq \nu_y(S(y,B(g y,r)))\leq C\cdot e^{-\beta\cdot\rho(y,g y)}.$$\elem 
Proof of Lemma \ref{SS} is by modification of Sullivan's argument in \cite[Proposition 3]{Sul}  in light of the appendix in \cite{BHM}; see Coornaert \cite[Proposition 6.1]{Coo}, Quint \cite[Lemma 4.10]{Qui}. 

We will use the following covering lemma.
\blem[$5r$-covering]\label{cover} Let $Y$ be a metric space. Given any collection of balls $\{B_\alpha\}_{\alpha\in\mathcal{J}}$, there exists a countable sub-collection of disjoint balls $\{B_i\}_{\alpha\in\mathcal{I}}$, $\mathcal{I}\subset\mathcal{J}$, such that $$\bigcup_{\alpha\in\mathcal{J}}B_\alpha\subset \bigcup_{i\in\mathcal{I}}5B_i,$$ where $5B_i$ denotes a ball with same centre as $B_i$ and five times the radius.
\elem

\blem[Shadow lemma for $\nu$]\label{shadow}
There exists $R_\nu>0$, such that for any $x \in X$, $r\geq R_\nu$, there exists $C_r\geq 1$, so that for all $y\in X$,
\beqnn\begin{split}\frac{1}{C_r}\cdot  e^{-\rho_G(x,y)} &\leq \nu_x(S_G(x,B_G(y,r)))  \leq C_r\cdot   e^{-\rho_G(x,y)}.\end{split}
\eeqnn  
\elem
\bproof
Set
\[k^H(y)=\begin{cases}\log g_H(a^{ \rho_X(y,\Gamma)}) &  n_y>0 \\ 0 & n_y=0,\end{cases} \]where $\rho(y,\Gamma):=\dist_\rho(y,\Gamma)$, for $\rho=\rho_X,\rho_G$.
 First assume that $x\in\Gamma$. We sketch the arguments in the situation when $y$ is in a horoball; if $y\in \Gamma$, the claim follows from Lemma \ref{SS} using \eqref{conformality}. Let $R_S>0$ be chosen according to Lemma \ref{SS}. Let $r\geq R_S+\delta_G$. Then, let $y\in H_g$. Let $y_1=gh_1\in gH$ be such that a geodesic ray $\gamma_y$ starting at $x$ intersecting a $\delta_X$-neighbourhood (by visuality of $(X,\rho_X)$) of $y$ enters $gH$ at $y_1$. Let $y_2=gh_2\in gH$ be the point where $\gamma_y$ exits the horoball $H_g$. We may assume that $$\min\,\{\rho_X(y_1,y),\rho_X(y,y_2)\}\geq 10\cdot r,$$ for Lemma \ref{SS}. 
\vspace{0.2cm}

\textbf{Step 1:} Write $$m_y=\min\,\{k\in\N\,\mid\, p_{xy'}\cap (gH\times \{k\})\neq \emptyset,\; \text{for all} \;y'\in B_G(y,r)\},$$ where $p_{xy'}$ is a preferred ray starting at $x$ passing through $y'$. First consider the case $m_y<\infty$.
In this case $$\max\,\{n_{y'}\mid\,y'\in B_G(y,r)\}< m_y+C_1,$$ where $C_1=C_1(X)>0$ is a constant larger than the maximum distance between geodesics and preferred rays.
Then there exists $r'=r'(r)>r$ such that any geodesic ray $\gamma_{y'}$ starting at $x$ and intersecting $B_G(y,r)$ at $y'$, intersects $U_H=U_{H}(y_2, a^{\rho_X(y,y_2)+r'})$, which is a ball of radius $a^{\rho_X(y,y_2)+r'}$ centred at $y_2$ in $gH$ as a subspace of  $Cay(\Gamma)$, with the word metric $d_{\Gamma}$ restricted to $gH$. Consider the collection of balls \beqnn  \quad \quad\hat{\mathcal{U}}_y=\{B_G(z,R_S+\delta_G) \,\mid\,z\in U_H, \; S_G(x,B_G(z,R_S+\delta_G))\subset S_G(x,U_{H})\},\eeqnn and apply the covering theorem Lemma \ref{cover} to obtain a disjoint subcover of balls with smaller radius $$\mathcal{U}_y:=\{B_G(z_i,R_S) \,\mid\,z_i\in U_{H}\}$$ such that $$\bigcup_i S_G(x,B_G(z_i,R_S))\,\subset\, S_G(x, U_H) \,\subset \,\bigcup_i S_G(x,B_G(z_i,10\cdot R_S)),$$ where the first inclusion is as a disjoint union. Moreover, $\#\, \mathcal{U}_y\approx_r g_H(a^{\rho_X(y,y_2)})$. 

Note that $$\nu_x(S_G(x,B_G(z_i, R_S)))\approx \nu_x(S_G(x,B_G(z_i,10\cdot R_S)))\approx e^{-\rho_G(x,y_2)},$$ by Lemma \ref{SS}.
Then, \beqn\label{case1}
\begin{split}
\nu_x(S_G(x, B_G(y,r)))& \approx  \sum_i \nu_x(S_G(x,B_G(z_i, R_S))) \\&
\approx_r g_H(a^{\rho_X(y,y_2)})\cdot \nu_x(S_G(x,B_G(z_1, R_S))) \\&
\approx_r g_H(a^{\rho_X(y,y_2)})\cdot e^{-\rho_G(x,y_2)} \\& 
\approx_r g_H(a^{\rho_X(y,y_2)})\cdot e^{-\rho_G(y,y_2)}\cdot e^{-\rho_G(x,y)} \\&
\approx_r e^{k^H(y)-\rho_G(y,\Gamma)}\cdot e^{-\rho_G(x,y)},
\end{split}
\eeqn
where we used Lemma \ref{triangle} and Theorem \ref{ancona}. 
\vspace{0.2cm}

\textbf{Step 2:} Next we consider the case $m_y=\infty$.
Note that in this case $|\rho_X(y,\Gamma)- \rho_X(y,y_1)|\leq C_2$, for some $C_2=C_2(\delta_X)>0$ (using the description of geodesics in a horoball). 

Note also that there is $\overline{r}=\overline{r}(r)$ such that $$S_G(x, B_G(y,r))\subset S_G(x, gH\setminus U_H(y_1,a^{\rho_X(y_1,y)-\overline{r}})),$$ and $$S_G(x, B_G(y,r))\supset S_G(x, gH\setminus U_H(y_1,a^{\rho_X(y_1,y)+\overline{r}})).$$
For $\overline{r}'\in [-\overline{r},\overline{r}]$, consider the decomposition $$gH\setminus U_H(y_1,a^{\rho_X(y_1,y)+\overline{r}'})= \bigcup_{i=0}^\infty \, A_H(y_1, i),$$ where 
$$A_H(y_1, i)=U_H(y_1,  a^{\rho_X(y_1,y)+\overline{r}'+(i+1)\cdot \overline{r}})\setminus U_H(y_1,  a^{\rho_X(y_1,y)+\overline{r}'+ i\cdot \overline{r}}).$$ 
For $i\in\N_0$, let $\mathcal{A}_i$ be a $\delta_G$-separated collection of balls $B(z^{(i)}_j,R_S)$ centred at points in $A_H(y_1,i)$, obtained as in the second step, such that $$\bigcup_{i=0}^\infty \bigcup_{\mathcal{A}_i} S_G(x, B_G(z^{(i)}_j,R_S)) \subset S_G(x,gH\setminus U_H(y_1,a^{\rho_X(y_1,y)+\overline{r}'})),$$ and $$S_G(x,gH\setminus U_H(y_1,a^{\rho_X(y_1,y)+\overline{r}'}))\subset \bigcup_{i=0}^\infty \bigcup_{\mathcal{A}_i} S_G(x, B_G(z^{(i)}_j,10\cdot R_S)).$$ 
Note that $$\#\mathcal{A}_i\approx_{\overline{r}} g_H(a^{\rho_X(y_1,y)+\overline{r}'+i\cdot \overline{r}})$$ ($\overline{r}$ may be chosen to be larger than an absolute constant to ensure this).

Consider the set of preferred paths for each $i\in\N_0$, $\mathcal{G}_i=\{p_{y_1z}: z\in A_H(y_1,i)\}$; they have lengths greater than $\dist_{\rho_X}(y_1, A_H(y_1,i))-c(\delta_X)$. Note that there is a geodesic ray $\hat{\gamma}$ starting at $x$ passing through $\delta_X$-neighbourhood of $y$ which hits the parabolic point $\xi_{gH}$ fixed by $gHg^{-1}$. Let $w_i$ be a nearest point projection to $\hat{\gamma}$ of a chosen point $w\in \bigcup \mathcal{G}_i$ (denoting the union of all sets in the collection $\mathcal{G}_i$; the union is not over $i$) such that $n_w$ is maximal in $\bigcup \mathcal{G}_i$. 

Then, using $\#\mathcal{A}_i\approx e^{k^H(w_i)}$, we note
\beqn\label{secondofs}
\begin{split}
\nu_x(S_G(x,B_G(y,r))) &\approx \sum_{i=0}^{\infty} \sum_{\mathcal{A}_i} \nu_x(S_G(x, B_G(z_j^{(i)},R_S))) \\ &
\approx \sum_{i=0}^{\infty} \sum_{\mathcal{A}_i} e^{-\rho_G(x,z^{(i)}_1)} \\&
\approx \sum_{i=0}^{\infty} \sum_{\mathcal{A}_i} e^{-\rho_G(x,y)}\cdot e^{-\rho_G(y,w_i)}\cdot e^{-\rho_G(w_i,z_1^{(i)})} \\&
\approx_r e^{-\rho_G(x,y)} \sum_{i=0}^{\infty} e^{k^H(w_i)-\rho_G(w_i,\Gamma)}\cdot e^{-\rho_G(y,w_i)} 
\end{split}
\eeqn
We used Lemma \ref{triangle} and Theorem \ref{ancona} (and the $\Gamma$-invariance of $\rho_G$) above. 

\vspace{0.2cm}

\textbf{Step 3:} Now we will show that $e^{k^H(y)-\rho_G(y,\Gamma)}\approx 1$. This simplifies \eqref{secondofs}. Note that $$\nu_x(S_G(x, B_G(y,r)))\approx_r e^{-\rho_G(x,y)}\cdot \nu_y(S_G(x, B_G(y,r)))\leq e^{-\rho_G(x,y)}.$$ The sources of shadows can be taken to be points in the visual boundary. Then by hyperbolicity of $X$, for the case addressed in Step 1, $$S_G(x,B_G(y,r-\delta_G))\subset S_G(\xi_{gH}, B_G(y,r))\subset S_G(x,B_G(y,r+\delta_G)),$$ where $\xi_{gH}\in\partial X$ is the fixed point of $gHg\inv$.
Then by comparison with \eqref{case1} we get 
\beqn\label{thirdofs}
1\gtrsim_r \nu_y(S_G(\xi_{gH}, B_G(y,r)))\approx_r e^{k^H(y)-\rho_G(y,\Gamma)},
\eeqn as $y$ is in a horoball and $r$ can be chosen arbitrarily. Note that \eqref{thirdofs} holds for any $y\in H_g$ as $x$ may be chosen suitably so that the case is as in Step 1.

 We may use \eqref{thirdofs} to get a similar lower bound as follows. We may assume without loss of generality $y=(e,n)$ for some $n\in\N$, $n\geq 10\cdot r$. 
Then note that for some $\overline{r}>0$, 
\beqn\label{fourthofs}
\nu_{(e,n)}(S_G(\xi_H,B_G((e,n),r))\approx_{\overline{r}} \nu_{(e,n)}(S_G(\xi_H, U_H(e,g_H(a^{\rho_X((e,n),e)+\overline{r}}))).
\eeqn
Note also that for any $m>n$, there is $\overline{r}>0$, such that 
\beqn\label{fifthofs}
\nu_{(e,n)}(S_G(\xi_H, H\setminus U_H(e,g_H(a^{\rho_X((e,n),e)+\overline{r}}))))\approx \nu_{(e,n)}(S_G((e,n),B_G((e,m),r))).
\eeqn
Then $\nu_{(e,n)}(S_G((e,n),B_G((e,m),r)))$ can be computed as in Step 2, by decomposing \newline$ H\setminus U_H(e,g_H(a^{\rho_X((e,m),e)+\overline{r}}))$ into suitable annuli $A_i$, and the hyperbolicity of $X$, doubling in $H$, as 
\beqn\label{sixthofs}
\begin{split}
\nu_{(e,n)}(S_G((e,n),B_G((e,m),r))) & \approx_r \sum_{i=0}^\infty e^{-\rho_G((e,n),(e,m+i\cdot\overline{r}))}\cdot \nu_{(e,m+i\cdot\overline{r})}(S_G((e,n),\bigcup \mathcal{A}_i)) \\&
\approx_r \sum_{i=0}^\infty e^{-\rho_G((e,n),(e,m+i\cdot\overline{r}))}\cdot e^{k^H((e,m+i\cdot\overline{r}))-\rho_G((e,m+i\cdot\overline{r}),e)}
\end{split}
\eeqn
Then, using the upper bound from \eqref{thirdofs}, given $c>0$ there is some constant $c'>0$, such that if $m\in N$ is smallest larger than $n+c'$, then $$\nu_{(e,n)}(S_G((e,n),B_G((e,m),r)))\leq c.$$ 

Thus by \eqref{fourthofs}, \eqref{fifthofs}, \eqref{sixthofs}, 
\beqn\label{seventhofs}
\nu_{(e,n)}(S_G(\xi_H,B_G((e,n),r))\gtrsim_r 1
\eeqn

Now we apply \eqref{thirdofs} and \eqref{seventhofs}, to \eqref{case1} and \eqref{secondofs} to get
$$\nu_x(S_G(x, B_G(y,r)))\approx_r e^{-\rho_G(x,y)}.$$
\vspace{0.2cm}

\textbf{Step 4:} Next assume $x\in X\setminus \Gamma$. Note that for any $R>0$ large enough, and $y\in X$ such that $\dist_{\rho_G}(x, B_G(y,2\cdot R))>20\cdot \delta_G$, there exists $x'\in\Gamma$ (in the coset associated to the horoball containing $x$), such that a geodesic ray starting at $x'$ satisfies $$S_G(x',B_G(y,R-5\cdot \delta_G))\subset S_G(x, B_G(y,R))\subset S_G(x',B_G(y,R+5\cdot \delta_G)),$$ and $(x'|y)_{x,\rho_G}\leq \tau,$ for an absolute constant $\tau>0$. Such a point $x'\in \Gamma$ exists by the description of geodesics in $X$ in \S \ref{construction}. Then the claim follows by harmonicity of $\nu$ and the estimate for shadows of balls seen from points in $\Gamma$.

This finishes the proof.
\eproof

\brem\label{imp_1}
Following are remarks on Lemma \ref{shadow}. Both points here will be used in coming sections.
\benum
\item From Step 4 of the proof (equations \ref{thirdofs} and \ref{seventhofs}), it follows that $e^{k^H((e,n))-\rho_G((e,n),e)}\approx 1$. Then if $d_H$ is the rank (of a finite index nilpotent-subgroup) of $H$, we get that $$G((e,n),e)\approx {(a^{d_H})}^{-n},\quad\lim_{n\to\infty}\,\frac{k^H((e,n))}{\rho_G((e,n),e)}= 1.$$ Then by Ancona's inequality and the description of geodesics in horoballs from \S \ref{construction}, we get that the Green's function inside a horoball $H$, up to absolute multiplicative constants depends only on the growth of $H$, as, $$\frac{G(x,y)}{m(y)}\approx {(a^{d_H})}^{-\rho_X(x,y)},$$ for all $x,y\in H_g$, for all $g\in\Gamma$, in other words $$|\rho_G(x,y)-d_H\cdot\log a\cdot\rho_X(x,y)|\leq C,$$ for some $C=C(X,P,m)>0$.
\item (The Patterson-Sullivan measure). Let $D_\Gamma$ be the critical exponent for $\Gamma$;
$$D_\Gamma:=\limsup_{n\to\infty}\,\frac{1}{n}\cdot \log(\#(\Gamma\cap B_X(e,n))).$$
 Then note that by Lemma \ref{QI}, $B_X(e,n)\subset B_G(e,C\cdot n)$, for some $C>0$. Thus $D_\Gamma<\infty$ by Lemma \ref{growth} below.
 
Let $\{\mu\}_x$ be a $D_\Gamma$-dimensional, $\Gamma$-quasiconformal Patterson-Sullivan density obtained by the Patterson construction.  The proof of Lemma \ref{shadow}  (Steps 1,2,4) can be applied also to $(\rho_X,\mu)$; we have for $x\in\Gamma$ and $y\in X$, and $r>R=R(\rho_X)>0$ (cf. \cite{SV}, \cite{BT}),\vspace{0.1cm}
$$\mu_x(S(x,B(y,r)))\approx_r e^{k_H(y)-D_\Gamma\cdot \rho_X(y,\Gamma)}\cdot e^{-D_\Gamma\cdot \rho_X(x,y)}.$$
For example, we check that \eqref{secondofs} translates to 
\beqnn
\begin{split}
\mu_x(S(x,B(y,r))) 
&\approx_r e^{-D_\Gamma\cdot\rho_X(x,y)} \sum_{i=0}^{\infty} e^{k_H(w_i)-D_\Gamma\cdot\rho_X(w_i,\Gamma)}\cdot e^{-D_\Gamma\cdot\rho_X(y,w_i)} \\ 
&\approx e^{-D_\Gamma\cdot\rho_X(x,y)}\cdot e^{k_H(y)-D_\Gamma\cdot\rho_X(y,\Gamma)}\sum_{i=0}^{\infty} e^{-(2\cdot D_\Gamma-d_H\cdot\log a)\cdot\rho_X(y,w_i)}.
\end{split}
\eeqnn
By finiteness of $\mu$, we get $$D_\Gamma>\frac{d_H}{2}\cdot \log a,$$ and thus $$\mu_x(S(x,B(y,r)))\approx_r e^{-D_\Gamma\cdot\rho_X(x,y)}\cdot e^{k_H(y)-D_\Gamma\cdot\rho_X(y,\Gamma)}.$$
By Remark \ref{imp_2} below, limit of supremums in defining $D_\Gamma$ may be replaced by the limit.
\eenum
\erem

\vspace{0.1cm}

\bproof[\textbf{Proof of Theorem \ref{char}}]
$(2)\implies (1)$ follows from the discussion in \S \ref{QC} and Theorem \ref{double}; the latter follows from Lemma \ref{shadow} and the (quasiruled) hyperbolicity of $(X,\rho_G)$. \eproof

We will need the asymptotic growth of $\Gamma$ in $X$ with respect to $\rho_G$.
\blem[Growth of $\Gamma$]\label{growth}
There exists $R_X>0$ such that, for any $x\in\Gamma$, and $2\leq n\in\N$,
\beqnn
\#\{\gamma\in\Gamma \,\mid\, n\cdot R_X\leq \rho_G(x,\gamma)\leq (n+1)\cdot R_X\} \approx e^{n\cdot R_X}.
\eeqnn
\elem

\bproof
Choose first $R_X>10\cdot R_\nu$. For each $n\in\N$, consider the set  $$\mathcal{V}_G(n)=\{g\in\Gamma \,\mid\, n\cdot R_X\leq \rho_G(x,g)\leq (n+1)\cdot R_X\},$$ and the collection $\{S_G(x,B_G(g, R_\nu))\}_{g\in \mathcal{V}_G(n)}$. Note that by hyperbolicity of $(X,\rho_G)$, $$\frac{1}{C}\sum_{g\in \mathcal{V}_G(n)}\ind_{S_G(x,B_G(g, R_\nu))}\leq 1,$$ where $C=C(R_X,\delta_G)>0$. Integrating with respect to $\nu_x$, we get 
\beqn\label{firstofgrowth}\sum_{g\in \mathcal{V}_G(n)} e^{-\rho_G(x,g)}\lesssim 1,\eeqn which is the required upper bound.

Suppose now that $g\in B_G(x,(n-1)\cdot R_X)$, for $2\leq n\in \N$ is such that 
$$\dist_{\rho_G}(x,gH)=\rho_G(x,g),$$ 
for some coset $gH$, that is, $g$ is a nearest point from $x$, of a horoball $H_{g}$ based at a parabolic fixed point, $\xi_{gH}$. Consider $H_{g'}\cap (B_G(x,(n+1)\cdot R_X)\setminus B_G(x, n\cdot R_X))$. Let $s_g$ be depth of the final horosphere of $H_{g'}$ intersecting $B_G(x,(n+1)\cdot R_X)\setminus B_G(x, n\cdot R_X)$. For the $j$-th horosphere of $H_{g}$, for $0\leq j \leq s_g$, consider the intersection 
$$A^g_j=(gH\times \{j\})\cap (B_G(x,(n+1)\cdot R_X)\setminus B_G(x, n\cdot R_X)).$$ 
Write $$U_H(n,j,c)=U_H\left(gH\times\{j\},\frac{n\cdot R_X -\rho_G(x,g)+j\cdot d_H\cdot\log a+c}{2\cdot d_H\cdot\log a}\right).$$
By the description of geodesics in $(X,\rho_X)$ inside horoballs (quasiruled-quasigeodesics in $(X,\rho_G)$) and Remark \ref{imp_1}, part one, we have 
$$A^g_j\subset U_H(n+1,j,c_1)\setminus U_H(n,j,c_2),$$ 
where $c_1, c_2>0$ are constants, and 
$$s_g= \frac{(n+1)\cdot R_X -\rho_G(x,g)}{d_H\cdot\log a}\pm c_3,$$ 
for a constant $c_3\geq 0$. Then by Lemma \ref{cover}, we have a collection of disjoint balls $\{B_G(z^j_i,R_\nu)\}_{1\leq i\leq i_j}$ of cardinality $i_j$ bounded by a constant times $\exp(\frac{n\cdot R_X-\rho_G(x,g)-j\cdot d_H\cdot\log a}{2})$, such that 
$$(gH\times \{j\})\cap (B_G(x,(n+1)\cdot R_X)\setminus B_G(x, n\cdot R_X))\subset \bigcup_{i=1}^{i_j} B_G(z^j_i,5\cdot R_\nu).$$ 
Consider the collection $$\mathcal{J}_g = \{B_G(z^j_i, R_\nu) \mid 0\leq j\leq s_g, 1\leq i\leq i_j\},$$ and note that 
\beqnn\begin{split}\#\mathcal{J}_g &\leq \sum_{j=0}^{s_g} e^{\frac{n\cdot R_X-\rho_G(x,g)-j\cdot d_H\cdot\log a}{2}}\\&
\lesssim  e^{\frac{n\cdot R_X-\rho_G(x,g)}{2}} \sum_{j=0}^{s_g} e^{\frac{-j\cdot d_H\log a}{2}} \\&
\lesssim e^{\frac{n\cdot R_X-\rho_G(x,g)}{2}}.
\end{split}\eeqnn  

Choose for $n\geq 1$, a set of coset representatives, 
$$\mathcal{T}_G(n)=\{g\in \mathcal{V}_G(n)\mid g\;\text{is a nearest point from}\;x\;\text{of a coset}\},$$ 
and for $n\geq 2$, $k_0\leq n-1$, 
$$\mathcal{W}_G(n,k_0)=\{g\in \mathcal{V}_G(n)\mid \text{coset of}\; H\; \text{containing}\;g\;\text{is outside}\; B_G(x,(n-k_0)\cdot R_X)\}.$$
Note that there exists a constant $C\geq 1$, such that 
\beqn\label{secondofgrowth}
1\leq \sum_{g\in \mathcal{W}_G(n,k_0)}\ind_{S_G(x,B_G(g, C\cdot R_\nu))}+\sum_{k=0}^{n-k_0}\sum_{g\in \mathcal{T}_G(k)}\sum_{\mathcal{J}_{g}}\ind_{S_G(x,B_G(z^j_i, C\cdot R_\nu))}
\eeqn
Integrating with respect to $\nu_x$ and using the upper bound \eqref{firstofgrowth},  
\beqnn\begin{split}
1\;\leq \;& C'\cdot \#\mathcal{W}_G(n,k_0)\cdot e^{-n\cdot R_X} \;+\;C''\cdot \sum_{k=0}^{n-k_0}\#\mathcal{V}_G(k)\cdot e^{\frac{n\cdot R_X-k\cdot R_X}{2}}\cdot e^{-n\cdot R_X}\\& \leq \; C'\cdot \#\mathcal{W}_G(n)\cdot e^{-n\cdot R_X} \;+\; C'''\cdot e^{\frac{-n\cdot R_X}{2}}\sum_{k=0}^{n-k_0} e^{\frac{k\cdot R_X}{2}},
\end{split}
\eeqnn
for constants $C', C'', C'''>0$. Then the claim follows by taking $k_0=2$, $R_X$ large such that $$C'''\cdot  e^{\frac{-n\cdot R_X}{2}} e^{\frac{(n-k_0+1)\cdot R_X}{2}}<1/2.$$
\eproof

Recall that for $n\geq 2$,$$\mathcal{V}_G(n)=\{g\in\Gamma \,\mid\, n\cdot R_X\leq \rho_G(x,g)\leq (n+1)\cdot R_X\},$$ and  $$\mathcal{W}_G(n,2)=\{g\in \mathcal{V}_G(n)\mid \text{coset of}\; H\; \text{containing}\;g\;\text{is outside}\; B_G(x,(n-2)\cdot R_X)\}.$$

The proof in fact gives the following two lemmas.
\blem[Growth of cosets]\label{growth_cosets}
There exists $R_X>0$ such that, for any $x\in\Gamma$, and $2\leq n\in\N$
\beqnn
\#W_G(n,2) \approx e^{n\cdot R_X}.
\eeqnn
\elem
\blem\label{mixing_boundary}
There exists $R_X>0$ such that, for any $x\in\Gamma$, $U$ any ball in $\partial X$, and $n\in\N$ large
\beqnn
\#\{g\in W_G(n,2) \,\mid\, S_G(x,B_G(\gamma,R_X))\cap U\neq\emptyset,\} \approx e^{n\cdot R_X}\cdot \nu_x(U).
\eeqnn
\elem
\bproof
The lemma follows by multiplying \eqref{firstofgrowth} and \eqref{secondofgrowth} with $\ind_U$, using the doubling property of $\nu_x$, and integrating with respect to $\nu_x$.
\eproof

\brem[Growth in $\rho_X$]\label{imp_2}
The proofs of Lemma \ref{growth}, Lemma \ref{growth_cosets} and Lemma \ref{mixing_boundary} apply also to $\rho_X$. The key ingredients are Sullivan's shadow lemma for a lower bound on measure of shadows of balls centred at orbit points, the second shadow lemma; Remark \ref{imp_1}, second part for upper bounds on the shadows of balls centred at non-orbit points, the polynomial growth of $H$ and hyperbolicity. The growth of $\Gamma$ and the growth of cosets for a ball $B_X(x,n\cdot R_X)$, for some $R_X>0$, in this case is of the order $e^{D_\Gamma\cdot n}$, and the ratio of the measure $\mu_x$ of the shadow of an annulus at distance $n\cdot R_X$ from $x$, with the measure of the union of the shadow of balls of radius $R_X$ centered at elements of $\Gamma$ which are coset representatives of horoballs in the annulus is uniformly bounded away from zero and infinity.
\erem

\section{The quotient walk}\label{dimensionharmonic}
We define the drift and entropy of the random walk in this section.
Given our reversible Markov chain $(X,P,m)$, we begin by considering the space $\widetilde{\Omega}:=(X^\N,\p_m)$, where $$\p_m=\sum_{z\in X} m(z)\cdot \p_z.$$ The measure $\p_m$ is invariant for the induced action of $\Gamma$ on $X^\N$. Then the quotient space $\Omega'=\Gamma\backslash X\times X^\N$ equipped with the quotient topology and $\sigma$-algebra, and measure $\p'_{\overline{m}}$, is a standard probability space, where $$\overline{m}:\{e\}\times \N_0\longrightarrow \{e\}\times \N_0,\quad \text{is defined by}\quad \overline{m}((e,k))=\frac{m((e,k))}{\sum_{n\in\N_0} m((e,n))},$$ and 
$$\p'_{\overline{m}}= \sum_{n\in\N_0}\overline{m}((e,n))\cdot \p_{(e,n)}.$$ 
The walk in $\Omega'$ can be interpreted as a (non-compact) topological Markov shift. \vspace{0.2cm}

Towards that end define $\pi_1: \Gamma\backslash X\times X^{\N}\longrightarrow X^{\N}$ by $$\{(e,n_0),(g_1,n_1),\ldots,(g_k,n_k),\ldots\}\overset{\pi_1}{\longmapsto} \{(g_1,n_1),\ldots,(g_k,n_k),\ldots\}.$$
and consider the quotient space (conditioning with respect to the zeroth-coordinate) $$(\Omega'', {(\pi_1)}_{\ast}\,\p'_{\overline{m}})=\pi_1\backslash(\Omega', \p'_{\overline{m}}).$$ One can define an (orbit jump) operation on $X$ that forgets the depth (orbit) of the first operand: $$(g_1,n_1)\ast (g_2,n_2)=(g_1g_2,n_2).$$ We define then the space of increments. Set 
\[ S_X = \left\lbrace (s,n)\in \Gamma\times\N_0\;\middle|\;
\begin{tabular}{@{}l@{}}
$s\in B_\Gamma(e,a)$, if $n=0$ and,\vspace{0.1cm}\\
$s\in B_H(e,a^{n+1})$, if $n\in\N$.
\end{tabular}
\right\rbrace \]

\brem[Definition of $S_X$] Here we picked a connected set of $X$ (containing $e$; exists by connectedness of the graph) as the set of representatives of $\Gamma\backslash X$, namely $\{e\}\times\N_0$. Then $$S_X=\{y\in X\,\mid\, p(x,y)\neq 0,\, x\in\Gamma\backslash X\}.$$
\erem
Note that there exists a map $\pi_2:S_X^{\N}\longrightarrow \Omega'',$ defined
\beqnn
\begin{split}
&\{(s_1^{(n_1)},n_1), (s_2^{(n_2)},n_2),\ldots,(s_k^{(n_k)},n_k),\ldots\} \overset{\pi_2}{\longmapsto} \\&
\{(s_1^{(n_1)},n_1), (s_1^{(n_1)},n_1)\ast(s_2^{(n_2)},n_2),\ldots,(s_1^{(n_1)},n_1)\ast\cdots\ast(s_{k}^{(n_{k})},n_{k}),\ldots\} \\&
=\left\lbrace(s_{1}^{(n_1)},n_1), (s_1^{(n_1)}s_2^{(n_2)},n_2), \ldots,(\prod_{l=1}^k s_{l}^{(n_l)}, n_l),\ldots\right\rbrace
\end{split}
\eeqnn
Note that $\pi_2$ is injective. This follows by observing that $X\xrightarrow{(g_1,n_1)\ast} X$ defined by \newline $(g_1,n_1)\ast (g_2,n_2)=(g_1g_2,n_2)$ is injective. Note also that $$\text{Supp}_{\Omega''}({(\pi_1)}_\ast\,\p'_{\overline{m}})\subset \pi_2(S_X^{\N}).$$ Then there is an isomorphism of measure spaces given by $\pi_2$: $$(\Omega'',{(\pi_1)}_\ast\,\p'_{\overline{m}})\overset{\pi_2\inv}{\longrightarrow} (S_X^{\N},{(\pi_2\inv)}_\ast\,{(\pi_1)}_\ast\,\p'_{\overline{m}})=:\Omega_1.$$
Write $$\p_{\overline{m}}={(\pi_2\inv)}_\ast\,{(\pi_1)}_\ast\,\p'_{\overline{m}}.$$ For a cylinder set $$C=[(s_1^{(n_1)},n_1),\ldots, (s_k^{(n_k)},n_k)]:=\{\omega\mid X_i(\omega)=(s_i,n_i), 1\leq i\leq k\}\subset S_X^{\N},$$ where $k\in\N$, by definition,
\beqnn
\begin{split}
\p_{\overline{m}}(C) &=\sum_{(e,n)\in\Gamma\backslash X} \overline{m}((e,n))\cdot p((e,n),(s_1^{(n_1)},n_1))\cdot \prod_{j=1}^{k-1}p\left(\left(\prod_{i=1}^{j}s_i^{(n_i)},n_i\right),\left(\prod_{i=1}^{j+1}s_i^{(n_i)},n_i\right)\right) \\& 
=\sum_{(e,n)\in\Gamma\backslash X} \overline{m}((e,n))\cdot p((e,n),(s_1^{(n_1)},n_1))\cdot \prod_{j=1}^{k-1} p((s_j^{(n_j)},n_j),(s_j^{(n_j)}s_{j+1}^{(n_{j+1})},n_{j+1})).
\end{split}\eeqnn
Define the left-shift $T$ on $\Omega_1$:
$$T(\{(s_1^{(n_1)},n_1),\ldots,(s_k^{(n_k)},n_k),\ldots\})=\{(s_2^{(n_2)},n_2),\ldots,(s_{k+1}^{(n_{k+1})},n_{k+1}),\ldots\}.\vspace{0.2cm}$$
Note that the increment coordinate-variables $\{X_i: S_X^{\N}\longrightarrow S_X\}_{i\in\N}$ are a Markov chain: \beqnn\begin{split}\p_{\overline{m}}\left[X_{n+1}=x_{n+1}\,\mid\, X_n=x_n,\ldots,X_1=x_1\right] &=\p_{\overline{m}}\left[X_{n+1}=x_{n+1}\,\mid\, X_n=x_n\right]\\& =p(x_n,x_n\ast x_{n+1}),\end{split}\eeqnn where $x_i\in S_X$.
We denote for $r\in\N$, $\prod_{i=1}^r X_i(\omega)=X_1(\omega)\ast\cdots\ast X_r(\omega)$.
\blem\label{ergodic_shift} 
$(\text{Supp}_{\Omega}(\p_{\overline{m}}), T, \p_{\overline{m}})$ is a mixing and invariant Markov shift.
\elem
\bproof
For the cylinder set $C=[(s_1^{(n_1)},n_1),\ldots,(s_k^{(n_k)},n_k)]$, 
$$T\inv(C)=\underset{(s,n)\in S_X}{\bigcup}\,[(s,n), (s_1^{(n_1)},n_1),\ldots,(s_k^{(n_k)},n_k)].$$ Then, writing $p(C)=\prod_{j=1}^{k-1} p((e,n_j),(s_{j+1}^{(n_{j+1})},n_{j+1}))$, we get
\beqnn
\begin{split} 
\p_{\overline{m}}(T\inv C)& = \sum_{(e,n')\in\Gamma\backslash X}\sum_{(s,n)\in S_X}\overline{m}((e,n'))\cdot p((e,n'),(s,n))\cdot p((e,n),(s_1^{(n_1)},n_1)) \cdot p(C)\\&
= \sum_{(e,n')\in\Gamma\backslash X}\sum_{(s,n)\in S_X} \overline{m}((e,n))\cdot p((e,n),(s\inv,n'))\cdot p((e,n),(s_1^{(n_1)},n_1)) \cdot p(C) \\&
= \sum_{(e,n)\in\Gamma\backslash X}\sum_{(s\inv,n')\in S_X} \overline{m}((e,n))\cdot p((e,n),(s\inv,n'))\cdot p((e,n),(s_1^{(n_1)},n_1)) \cdot p(C) \\&
=\sum_{(e,n)\in\Gamma\backslash X} \overline{m}((e,n)) \cdot p((e,n),(s_1^{(n_1)},n_1))\cdot p(C)= \p_{\overline{m}}(C).
\end{split}
\eeqnn

Mixing follows since the Markov shift is topologically mixing, $\p_{\overline{m}}$ is Markov and $$p_{(s,n)}:= \sum_{(e,n')\in \Gamma\backslash X}\, \overline{m}((e,n'))\cdot p((e,n'),(s,n))$$ defines a stationary initial probability distribution on $S_X$.

\eproof

\blem\label{Kingman1} For $\p_{\overline{m}}$-a.e $\omega\in S_X^{\N}$, the following limits exist:
$$\underset{n\to\infty}{\lim}  \frac{\rho_X(X_1(\omega),\prod_{i=1}^{n+1} X_i(\omega))}{n}=l,$$
$$\underset{n\to\infty}{\lim}  \frac{\rho_G(X_1(\omega),\prod_{i=1}^{n+1} X_i(\omega))}{n}= l_G.$$
\elem
\bproof
We have \beqnn
\begin{split}
\rho_G(X_1(\omega),\prod_{i=1}^{r+s+1} X_i(\omega))&
\leq \rho_G(X_1(\omega),\prod_{i=1}^{r+1} X_i(\omega))+\rho_G(\prod_{i=1}^{r+1} X_i(\omega),\prod_{i=1}^{r+s+1} X_i(\omega)) \\&
= \rho_G(X_1(\omega),\prod_{i=1}^{r+1} X_i(\omega))+\rho_G(X_1(T^r\omega),\prod_{i=1}^{s+1} X_i(T_r\omega)),
\end{split}
\eeqnn
where the equality is by the definition of $\ast$ and the $\Gamma$-invariance of $\rho_G$.
Moreover, writing $$l_G(r)(\omega):=\frac{\rho_G(X_1(\omega),\prod_{i=1}^{r+1} X_i(\omega))}{r},$$ we have by the Markov property (and \S \ref{green_metric}; Harnack inequality \ref{harnack2}),
\beqnn
\begin{split}
\int_{\Omega_1} l_G(1)(\omega) \,d\p_{\overline{m}}(\omega)&
=\sum_{x_1\in S_X}\mathbb{E}_{\p_{\overline{m}}}\left[\ind_{\{X_1=x_1\}}\cdot \rho_G(X_1,X_1\ast X_2)\right] \\&
=\sum_{x_1\in S_X} \sum_{z_2\in X}\rho_G(x_1,z_2) \cdot p(x_1,z_2) \cdot \p_{\overline{m}}[X_1=x_1] \approx 1.
\end{split}
\eeqnn
Similar holds for $\rho_X$. The lemma follows by Lemma \ref{ergodic_shift} and the Kingman subadditive ergodic theorem.
\eproof

\blem\label{Kingman_2}
For $\p_{\overline{m}}$-a.e $\omega\in S_X^{\N}$, the following limit exists:
$$
\lim_{n\to\infty} \frac{-\log (p^{(n)}(X_1(\omega),\prod_{i=1}^{n+1}X_i(\omega)))}{n}=:h
$$
\elem
\bproof
Note that $h(n,P):= -\log (p^{(n)}(X_1(\omega),\prod_{i=1}^{n+1}X_i(\omega)))$ is subadditive. Also, 
\beqnn
\begin{split}
\int_{\Omega_1} h(1,P)(\omega)\,d\p_{\overline{m}}(\omega) &
=-\int_{\Omega_1} \log p^{(1)}(X_1(\omega),X_1(\omega)\ast X_2(\omega))\,d\p_{\overline{m}}(\omega) \\&
=\sum_{x_1\in S_X} \sum_{z_2\in X} -\log p(x_1,z_2)\cdot p(x_1,z_2)\cdot \p_{\overline{m}}[X_1=x_1] \\&
=\sum_{x_1\in S_X} \sum_{z_2\in X} -\log m(x_1)\cdot p(x_1,z_2)\cdot \p_{\overline{m}}[X_1=x_1] \,+ O(1)\\&
=\sum_{n\in\N_0}\sum_{x_1\in S_X} -\log m(x_1)\cdot \overline{m}((e,n))\cdot p((e,n),x_1) \,+ O(1)\\&
=\sum_{n\in\N_0}-\log m((e,n))\cdot \overline{m}((e,n))\sum_{x_1}p((e,n),x_1) \,+O(1)\\&
=\sum_{n\in\N_0}-\log m((e,n))\cdot \overline{m}((e,n))\,+O(1)<\infty.
\end{split}
\eeqnn
Then the claim follows again by Kingman subadditive ergodic theorem and Lemma \ref{ergodic_shift}.
\eproof

\vspace{0.2cm}
Let $\{Z_i:\Omega''\longrightarrow X\}_{i\in\N}$, be the coordinate variables in $\Omega''$. Then $Z_i\circ\pi_2=\prod_{i=1}^r X_i$.

\bproof[\textbf{Proof of Theorem \ref{drift}}]
The corresponding limits exist for the quotient walk. Indeed, for $(\Omega'',{(\pi_1)}_\ast\,\p'_{\overline{m}})$ we have by Lemma \ref{Kingman1}
$$ \underset{n\to\infty}{\lim}  \frac{\rho_X(Z_1,Z_{n+1})}{n}=l,\quad \underset{n\to\infty}{\lim}  \frac{\rho_G(Z_1,Z_{n+1})}{n}= l_G,$$ and by Lemma \ref{Kingman_2}
$$\underset{n\to\infty}{\lim} \frac{-\log (p^{(n)}(Z_1,Z_{n+1}))}{n}=:h,$$ $\p'_{\overline{m}}$-a.e in $\Omega''$. Consider the topological inclusion $\Omega'\subset \Omega''$ and note that $\p'_{\overline{m}}\ll {\pi_1}_\ast\p'_{\overline{m}}\restr \Omega'$, where ${\pi_1}_\ast\p'_{\overline{m}}\restr \Omega'$ is the restriction of ${\pi_1}_\ast\p'_{\overline{m}}$ to the subspace $\Omega'$. Thus the limits exist a.e in $(\Omega',\p'_{\overline{m}})$.

Returning to the Markov chain $(X,P,m)$, we have that $(X^\N,\p_m)$ is the disjoint union of the translates (by coset representatives of $H$) of the preimage of $(\Omega',\p'_{\overline{m}})$. Then the theorem follows from $\Gamma$-invariance of $(X,P,m)$.
\eproof
\vspace{0.2cm}

\underline{Ergodicity of $(\Gamma,\nu\otimes\nu$)}: Consider the two-sided `shift' space $(\Omega_2:=Supp_{S_X^{\Z}}(\tilde{\p}_{\overline{m}}),T,\tilde{\p}_{\overline{m}})$, where the probability measure $\tilde{\p}_{\overline{m}}$ is defined by its values on cylinder sets 
\beqn
\begin{split}
C& =[(s_{-i}^{-1},n_{-i}),\ldots,(s_{-1}^{-1},n_{-1}),(e,n_0),(s_1,n_1),\ldots,(s_j,n_j)] \\ &:=\{\omega\in S_X^{\Z}: X_{-l_1}(\omega)=(s_{-l_1}^{-1},n_{-l_1}), X_0(\omega)=(e,n_0), X_{l_2}(\omega)=(s_{l_2},n_{l_2}), 1\leq l_1\leq i, 1\leq l_2\leq j\},\end{split}\eeqn
 where $i,j\in\Z$ and $X_l$ are the coordinate-random variables,
$$\tilde{\p}_{\overline{m}}(C)=\prod_{l_1=0}^{i-1}\, p((e,n_{-l_1}),(s_{-(l_1+1)}^{-1},n_{-(l_1+1)})) \cdot \overline{m}((e,n_0))\cdot \prod_{l_2=0}^{j-1}\,p((e,n_{l_2}),(s_{l_2+1},n_{l_2+1})).$$
The map $T$, defined for $$x=\{\ldots, (s_{-i}^{-1},n_{-i}),\ldots (s_{-1}^{-1},n_{-1}),(e,n_0),(s_1,n_1),\ldots, (s_j,n_j),\ldots\},$$ by
$$T(x)=\{\ldots, (s_{-i}^{-1},n_{-i}),\ldots (s_{-1}^{-1},n_{-1}), (s_1^{-1},n_0),(e,n_1), (s_2,n_2),\ldots (s_j,n_j),\ldots \},$$ is a bijection invariant for $\tilde{\p}_{\overline{m}}$. Indeed, for the set $C$ above,
$$T\inv C=\bigcup_{(s\inv,n)}[(s\inv,n),(s_{-i}\inv,n_{-i})\ldots, (s_{-2}\inv,n_{-2}),(e,n_{-1}), (s_{-1},n_0), (s_1,n_1) \ldots, (s_j,n_j)],$$ and 
\beqn
\begin{split}
\tilde{\p}_{\overline{m}}(T\inv C) & = \sum_{(s^{-1},n)} p((e,n_{-i}),(s^{-1},n))\cdot \prod_{l_1=1}^{i-1}\, p((e,n_{-l_1}),(s_{-(l_1+1)}^{-1},n_{-(l_1+1)})) \\&
 \times \overline{m}((e,n_{-1})) \cdot  p((e,n_{-1}),(s_{-1},n_0))\cdot \prod_{l_2=0}^{j-1}\, p((e,n_{l_2}),(s_{l_2+1},n_{l_2+1})) \\ &
= \tilde{\p}_{\overline{m}}(C),
\end{split}
\eeqn
using the reversibility $(P,m)$ and $\Gamma$-invariance. For the random variables $Z_n=X_0\ast X_1\ast\cdots\ast X_{n-1}$,  $Z_{-n}=X_0\ast X_{-1}\ast\cdots \ast X_{-(n-1)}$ and $Z_\infty=\lim_{n\to\infty} Z_{n}$, $Z_{-\infty}=\lim_{n\to\infty} Z_{-n}$ define the map
$$\pi_3: (\Omega_2,\tilde{\p}_{\overline{m}})\longrightarrow (\partial X\times \partial X, m\cdot\nu\otimes\nu),$$ where 
$$m\cdot\nu\otimes\nu:=\sum_{n\in \N_0} \overline{m}((e,n))\cdot \nu_{(e,n)}\otimes \nu_{(e,n)},$$ given by 
$$\pi_3(\omega)=(Z_{-\infty}(\omega),Z_\infty(\omega)).$$
Note that $\pi_3$ is measure-preserving. Then observing that for a $\Gamma$-invariant $m\cdot\nu\otimes\nu$-measurable set $A$ of positive measure, $\pi_3\inv(A)$ is a $T$-invariant $\tilde{\p}_{\overline{m}}$-measurable set of positive measure; and by ergodicity of $(\Omega_2,T,\tilde{\p}_{\overline{m}})$, we note the following (see \cite[Theorem 6.3]{Kai2} for the case of walks in Cayley graphs).
\blem The action $(\partial X\times \partial X, \Gamma, m\cdot\nu\otimes\nu)$, given by $g(\xi,\eta)=(g\xi,g\eta)$ for $(\xi,\eta)\in\partial X\times \partial X$ is ergodic.
\elem

\subsection{Cusp excursions and geodesic tracking}

In this section we write a geodesic tracking statement and a cusp-excursion statement which we will be needed in the next sections.

Apply the measurable selection theorem (and Arzela-Ascoli) to associate measurably to each $\xi\in\partial X$, a $\rho_X$-geodesic ray $\gamma^{x,\xi}$ starting at $x$, such that $\gamma^{x,\xi}_\infty=\xi$. The proof is a modification of an argument in \cite[Proposition 3.8]{BHM} for independent increment random walks in hyperbolic groups.

\blem\label{tracking}
Let $x\in X$. Let $k\in\N$. Then for any $D>0$, $$\p_x[\rho_G(Y_k,\gamma^{x,Y_\infty})\geq D]\leq C\cdot e^{- D},$$ for some $C>0$. In particular, $\p_x$-a.e $\omega\in X^{\N}$,
$$\limsup_{n\to\infty} \frac{\rho_G(Y_n(\omega), \gamma^{x,Y_\infty(\omega)})}{\log n}\leq 1.$$
\elem
\bproof
We get by the Markov property,
\beqnn
\begin{split}
\p_x[\rho_G(Y_k,\gamma^{x,Y_\infty})\geq D] &
=\sum_{y\in X} \p_{x}[Y_k=y, \rho_G(y,\gamma^{x,Y_\infty})\geq D] \\&
\leq \sum_{y\in X} \p_{x}[Y_k=y,Y_\infty\in S_G(y,B_G(x_{y},R))] \\&
= \sum_{y\in X} \p_{x}[Y_k=y]\cdot \p_x[Y_\infty\in S_G(y,B_G(x_{y},R))\,|\, Y_k=y] \\&
= \sum_{y\in X} \p_{x}[Y_k=y]\cdot \nu_{y}(S_G(y,B_G(x_{y},R))), 
\end{split}
\eeqnn
where $x_{y}$ is a point on a $\rho_X$-geodesic $[x,Y_k]$ at $\rho_G$-distance $D-\delta_G$ from $Y_k$, and the second inequality above is from the triangle inequality applied to the (formed of measurably chosen $\rho_X$-geodesics) $\rho_G$-quasigeodesic triangle $[x,Y_k, Y_\infty]$, once $R>0$ is chosen suitably. Then
$$\p_x[\rho_G(Y_k,\gamma^{x,Y_\infty})\geq D] \leq \sum_{y\in X} \p_{x}[Y_k=y]\cdot e^{-\rho_G(y,x_y)}\leq C\cdot e^{-D},$$ for some $C>0$. 

Next, note in particular that for $n\in\N$ and $\epsilon>0$, $$\p_x[\rho_G(Y_n,\gamma^{x,Y_\infty})\geq (1+\epsilon)\cdot\log n]\leq C\cdot \frac{1}{n^{1+\epsilon}},$$ so the second part of the claim follows from the first Borel-Cantelli lemma.
\eproof
\brem\label{tracking2}
Lemma \ref{tracking} implies by Lemma \ref{QI}, analogous statements for $\rho_X$, that is, there exist $b,C>0$ such that for any $k\in\N$, $D>0$, $$\p_x[\rho_X(Y_k,\gamma^{x,Y_\infty})\geq D]\leq C\cdot e^{- b\cdot D}\quad\text{and}\quad \limsup_{n\to\infty} \frac{\rho_X(Y_n(\omega), \gamma^{x,Y_\infty(\omega)})}{\log n}<\infty,$$ for $\p_x$-a.e $\omega\in X^{\N}$.
\erem

Next is a cusp excursion result.
 \blem[Cusp-excursion]\label{cusp_excur} 
We have for $x\in X$,
$$\nu_x\left(\left\lbrace\gamma_{\infty}^x\in \partial X\;\middle|\; \underset{t\to\infty}{\limsup} \,\frac{\rho_G(\gamma^x(t),\Gamma)}{\log t}=1\right\rbrace\right)=1,$$ and 
$$\mu_x\left(\left\lbrace\gamma_{\infty}^x\in \partial X\;\middle|\; \underset{t\to\infty}{\limsup} \,\frac{\rho_X(\gamma^x(t),\Gamma)}{\log t}=\frac{1}{2\cdot D_\Gamma-d_H\cdot \log a}\right\rbrace\right)=1.$$
where $\gamma_x$ is any $(X,\rho_X)$-geodesic ray starting at $X$ (the first limit is on times where $\gamma^x(t)\in V_\Gamma$).
\elem
\bproof
First suppose $x\in\Gamma$.
Let us denote by $S_H$ a chosen set of coset representatives of $H$, such that for $g\in S_H$, $$\rho_{G}(x,g)=\dist_{\rho_G}(x,H_g).$$ Write $l(g)=\rho_G(x,g)$ for $g\in S_H$. Let $f:[2,\infty)\to(0,\infty)$ be an increasing Lipschitz function. We think of the function $f$ as prescribing the maximum depth that a geodesic ray starting at $x$ reaches in a horoball in terms of the distance of the horoball to $x$. For each $g\in S_H$, choose a point $p_f(g)\in H_g$, such that $p_f(g)$ lies on a geodesic joining $x$ to $\xi_{gH}$, via a $\delta_G$-neighbourhood of $g$, such that $|\rho_G(g,p_f(g))-f(l(g))|\leq C$, where $C=\max \{\rho_G(x,y):\rho_X(x,y)\leq 1\}$.
Consider the set $$B_f=\limsup\,\{S_G(x,B_G(p_f(g),R))\mid\,g\in S_H\}.$$ Next recall the collections $\mathcal{A}_G(n)$ from Lemma \ref{growth_cosets} (the representatives chosen to be elements of $S_H$). For any $n\in\N$ and $m\in\N$ large (depending on $n$ and $f$) note that by Lemma \ref{mixing_boundary} for $g\in \mathcal{A}_G(n)$,
\beqn\label{firstofcusp_excur}
\#\{g'\in \mathcal{A}_G(m)\mid\, S_G(x,B_G(g', R))\cap S_G(x, B_G(p_f(g),R))\neq \emptyset\}
\lesssim e^m\cdot \nu_X(S_G(x, B_G(p_f(g),R))).
\eeqn
Write $A^f_G(n)=\bigcup\,\{S_G(x,B_G(p_f(g),R))\,\mid\, g\in \mathcal{A}_G(n)\}$. Then we have by \eqref{firstofcusp_excur},
$$\nu_x(A^f_G(n)\cap A^f_G(m))\lesssim \nu_x(A^f_G(n))\cdot \nu_x(A^f_G(m)).$$ 
Then by a Borel-Cantelli lemma, we have $\nu_x(B_f)>0$ if and only if $$\sum_{i\in N}\nu_x(A^f_G(n_i))=\infty,$$ for some sequence $\{n_i\}_i$, for which $\{n_i,n_{i+1}\}$, for all $i\in\N$ large, satisfy \eqref{firstofcusp_excur} (with $n=n_i$ and $m=n_{i+1}$). Moreover, by hyperbolicity of $(X,\rho_G)$, the ergodicity of the $\Gamma$ action on $(\partial X,\nu)$ (see Remark \ref{ergodicity}), and by an argument which demonstrates a $\Gamma$-invariant subset of positive measure contained in $B_f$ (for example \cite[p. 821]{HP}), we have $$\nu_x(B_f)=1,\quad \text{if}\quad \sum_{n\in \N}\nu_x(A^f_G(n))=\infty;\quad\nu_x(B_f)=0,\quad \text{if}\quad \sum_{n\in \N}\nu_x(A^f_G(n))<\infty.$$

Next note (using Lemma \ref{shadow}, Lemma \ref{growth}) that for the function 
$$f(t)= \log t,\quad\text{we have}\quad\sum_{i\in \N}\nu_x(A^f_G(n_i))\approx \sum_i \frac{1}{n_i}=\infty,$$ 
for a sequence $\{n_i\}$ which satisfies quasi-independence and some $\tilde{n}\in\N$ depending on $f$, and for any $\epsilon>0$, and 
$$f_\epsilon(t)= (1+\epsilon)\cdot \log t,\quad\text{we have}\quad \sum_{n\in \N}\nu_x(A^{f_\epsilon}_G(n))\lesssim \sum_n \frac{1}{n^{1+\epsilon}}<\infty.$$
For these choices of $f$ and $f_\epsilon$, $\nu_x(B_f\setminus B_{f_\epsilon})=1$. For $\xi\in B_f\setminus B_{f_\epsilon}$, and any $\rho_X$-geodesic $\gamma^x$, starting at $x$, with $\gamma^x(\infty)=\xi$, there is an increasing sequence $t_n$, $n\in\N$ of times when the geodesic $\gamma^x$ reaches far inside the horoballs $H_{g_n}$ such that for $n$ large, $t_n\in (k_n+f(k_n)-c_1,k_n+f_\epsilon(k_n)+c_2)$, for $k_n\in\N$, where $|l(g_n)-k_n|\leq 1/2$ and $c_1,c_2>0$ are constants.
Then since $\gamma^x$ is quasiruled,
\beqnn
\log(k_n)-c_1\leq  \dist_{\rho_G}(\gamma^x(t_n),\Gamma)\leq  (1+\epsilon) \cdot \log(k_n)+c_2,
\eeqnn 
for all $\epsilon>0$. It follows by varying $\epsilon$ that (see Remark \ref{imp_1}) $$\limsup_{t\to\infty} \frac{\rho_G(\gamma^x(t),\Gamma)}{\log t}=1.$$

The proof for $\rho_X$ is similar. In this case the functions are 
$$f^X(t)= \frac{\log t }{2\cdot D_\Gamma-d_H\cdot \log a}, \quad\text{and}\quad f_\epsilon^X(t)= \frac{(1+\epsilon)\cdot \log t}{2\cdot D_\Gamma-d_H\cdot \log a},$$ for $\epsilon>0$ and $\rho_G$ is replaced by $\rho_X$.

For general $(e,k)$, $k\in\N$, $$\nu_{(e,k)}(S_G((e,k),B(y,R_X)))\approx_k e^{-\rho_G(e,y)} \quad \text{and} \quad\rho_G((e,k),y)=\rho_G(e,y)+O(1),$$ for all $y\in X$. An identical argument then gives the same logarithmic laws for depths of geodesics starting at arbitrary $x\in X$.
\eproof

Before proving Theorem \ref{dimdriftent}, in the next section we digress to compare the asymptotic quantities defined for the random walk. 

\section{Entropy and drift}\label{guiver} By Remarks \ref{imp_1} we have $0<D_\Gamma<\infty$. Note that $h$ is non-zero by Lemma \ref{8to1}, Lemma \ref{basic_compar} and finite since $h\leq \mathbb{E}_{\p_{\overline{m}}}[h(1,P)]<\infty$). From Theorem \ref{dimdriftent}, $h=l_G$.

\brem Let $\lambda$ be the supremum of all constants for which the estimate of Lemma \ref{basicQI} holds for all pairs in $X$. 
We have $$-\log \rho(P)\;\leq\; \min\left\lbrace\frac{h - l\cdot\log \lambda}{1+l}, -\log \lambda\right\rbrace.$$ 

Indeed, note that for any $\epsilon>0$, there is $n=n_\epsilon\in\N$ large and $\omega\in X^{\N_0}$ (with $Y_0(\omega)=x$) such that 
\beqn
\begin{split}
h+\epsilon \geq \frac{-\log p^{(n)}(x,Y_n(\omega))}{n} &\geq -\frac{1}{n}\log\left(\frac{1}{m(x)}\right)-\log\left(\frac{1}{\lambda}\right)\cdot\frac{\rho_X(x,Y_n(\omega))}{n}+\frac{-\log p^{(n+r(x,Y_n(\omega)))}(x,x)}{n} \\&
\geq \epsilon -\log\left(\frac{1}{\lambda}\right)\cdot (l+\epsilon) + (-\log \rho(P)-\epsilon)(1+l-\epsilon).
\end{split}
\eeqn
It follows from Lemma \ref{basicQI} and Corollary \ref{8to1}, that $-\log \rho(P)\leq -\log \lambda$.
\erem

\bprop\label{fund} We have $$l/C\,\leq\, h\,\leq\, l\cdot D_\Gamma$$ for $(X,P,m)$, where $C=C(P,m)$ is a constant.\eprop
\bproof
Let $x\in\Gamma$. By Theorem \ref{drift}, Remark \ref{tracking2}, Lemma \ref{cusp_excur} and Egorov's theorem, given $\epsilon>0$, there is $N_\epsilon>0$ and a measurable set $\Omega_\epsilon\subset X^{\N}$ with $\p_x(X^{\N}\setminus \Omega_\epsilon)<1/2$, such that for $\omega\in \Omega_\epsilon$, 
\benum
\item $\rho_X(x,Y_n(\omega))\in B_X(x,(l+\epsilon)\cdot n)\setminus B_X(x,(l-\epsilon)\cdot n),$
\vspace{0.2cm}
\item $p^{(n)}(x,Y_{n}(\omega))\leq e^{-(h-\epsilon)\cdot n},$
\vspace{0.2cm}
\item $\rho_G(\gamma^{x,Y_\infty(\omega)}(n),\Gamma)\leq C\cdot \log n,$
\vspace{0.2cm}
\item $\rho_X(Y_n(\omega),\gamma^{x,Y_\infty(\omega)})\leq C\cdot \log n,$
\eenum
for any $n>N_\epsilon$ and a constant $C=C(X,P,m)>0$. For each $g\in \Gamma$ and $k\in\N$, let $\{B_i^k\}_i$ be a collection of $m$-nodes which cover $gH\times\{k\}$ and any point lies in at most $c$ sets of the form $B_i^k$, for some constant $c\geq 1$ (got by applying Lemma \ref{cover} to a suitable collection of balls in $H$). Let $B_i^0=g_i$ be singletons from $\Gamma$. Denote by $y_i^k$ the centre of $B_i^k$. Write $$A_n=(B_X(x,(l+\epsilon)\cdot n)\setminus B_X(x,(l-\epsilon)\cdot n))\bigcap\{y\in X\mid \exists \omega\in\Omega_\epsilon, Y_n(\omega)=y\}.$$ 
First note that if $B_i^k\cap A_n\neq\emptyset$, then for some $\omega\in \Omega_\epsilon$, $$m(y_i^k)=m(Y_n(\omega))=\frac{1}{e^{\log a\cdot d_H\cdot \rho_X(Y_n(\omega),\Gamma)}},$$ and by choice of $\Omega_\epsilon$, that is, using (3) and (4) and Lemma \ref{QI}, 
$$\rho_X(Y_n(\omega),\Gamma))\leq  C\cdot\log ((l+\epsilon)\cdot n),$$ for some constant $C=C(X,P,m)>0$. So
$$m(y_i^k) \geq \frac{1}{((l+\epsilon)\cdot n)^c},$$ for some constant $c=c(X,P,m)>0$, and by choosing $n$ large enough $$m(y_i^k)\geq e^{-n\cdot\epsilon}.$$
Next note that
\beqnn
\begin{split}
c/2\leq c\cdot \sum_{y\in A_n}p^{(n)}(x,y) &\leq \sum_k \sum_{B_i^k\cap A_n\neq\emptyset} \sum_{y\in B_i^k\cap A_n}p^{(n)}(x,y) \\ &
\leq \sum_k \sum_{B_i^k\cap A_n\neq\emptyset} e^{-n\cdot(h-\epsilon)}\cdot \#B_i^k  \\ & 
\lesssim \sum_k \sum_{B_i^k\cap A_n\neq\emptyset} e^{-n\cdot(h-\epsilon)} \cdot \frac{1}{m(y_i^k)} \leq  \sum_k \sum_{B_i^k\cap A_n\neq\emptyset} e^{-n\cdot(h-2\epsilon)}.
\end{split}
\eeqnn 
Also by second part of Remark \ref{imp_1}, 
$$\sum_k \sum_{B_i^k\cap A_n\neq\emptyset} e^{(d_H\cdot\log a-D_\Gamma)\cdot\rho_X(y_i^k,\Gamma)}\cdot e^{-D_\Gamma\cdot\rho_X(x,y_i^k)}\lesssim 1,$$ and 
$$e^{-D_\Gamma\cdot\rho_X(x,y_i^k)}\geq e^{-D_\Gamma\cdot (l+\epsilon)\cdot n}, \quad e^{(d_H\cdot\log a-D_\Gamma)\cdot\rho_X(y_i^k,\Gamma)}\geq e^{-D_\Gamma\cdot\rho_X(y_i^k,\Gamma)}\geq e^{-D_\Gamma\cdot n\cdot\epsilon},$$ where in the last lower bound, we used that $\rho_X(Y_n(\omega),\Gamma)$ is at most $O(\log n)$. Thus we get 
$$\sum_k \sum_{B_i^k\cap A_n\neq\emptyset} 1\,\leq\, e^{D_\Gamma\cdot(l+2\epsilon)\cdot n}.$$ Then

$$e^{(h-2\epsilon)\cdot n} \lesssim e^{D_\Gamma\cdot (l+2\epsilon)\cdot n},$$ which leads to the upper bound. 

For the lower bound note that for any $\epsilon>0$, on a $\p_e$-positive measure set for large $n$, $$(l-\epsilon)\cdot n\leq \rho_X(e,Y_n)\leq C\cdot \rho_G(e,Y_n)\leq C\cdot (l_G+\epsilon)\cdot n,$$ by Lemma \ref{QI}, for some $C>0$.
\eproof

The characterization of the situation the equality $h=l\cdot D_\Gamma$ holds, in terms of the Identity map being a $(1,C)$-quasiisometry between the metrics $D_\Gamma\cdot\rho_X$ and $\rho_G$ and equivalence of the harmonic and Patterson-Sullivan measures $\nu$ and $\mu$ as in \cite[Theorem 1.5]{BHM}, is also true in the present context. This uses the shadow lemmas for $\mu$ and $\nu$ and the ergodicity of the action of $\Gamma$ on $\partial X\times\partial X$. The proof is by suitable modification of the arguments of \cite{BHM} given the existence of the drifts $l_G$, $l$, entropy $h$. 

Some obstructions to equality are then clear from the shadow lemmas. Let the rank of a finitely-generated virtually nilpotent subgroup mean the rank of a finite-index nilpotent subgroup, we observe the following.
\bcor\label{rank}
Let $\Gamma$ be a finitely-generated group, hyperbolic relative to finitely many, finitely-generated, infinite, virtually nilpotent subgroups. Let $X_\Gamma$ be the corresponding cusped graph. Then there is a random walk $(X,P)$ such that, if $h=l\cdot D_\Gamma$ for $(X,P)$, then all the parabolic subgroups have rank $D_\Gamma/\log a$, where $D_\Gamma$ is the critical exponent for the action of $\Gamma$ in $X$ and $a>1$ is an absolute constant.
\ecor
\bproof[Proof of Corollary \ref{rank}] Consider the space $X=X_\Gamma(e)$. 
Recall from Remark \ref{imp_1} that, $$|\rho_G(x,y)-d_H\cdot\log a\cdot\rho_X(x,y)|\leq C,$$ for $x,y$ in horoballs and some constant $C>0$. Thus $h=l\cdot D_\Gamma$ would imply that all the parabolic subgroups have the same rank $d_H$ and $D_\Gamma=d_H\cdot \log a$, that is the rank of the parabolic subgroups is the critical exponent of the action of $\Gamma$ in $X$ (by the shadow lemmas; Lemma \ref{shadow}, part two of Remark \ref{imp_1}). \eproof Note that this is case the Patterson-Sullivan measure is Ahlfors-Regular. 
\vspace{0.3cm}

Next recall that $\Omega_1$ is the space of sequences of increments for the quotient walk with probability measure $\p_{\overline{m}}$; the walk denoted $X_i$. Recall that the depth function (which specifies the orbit) of a point $x\in X$, is denoted $n_x$.
\bprop\label{intrep} We have
\beqnn
\begin{split}
h=l_G &=\int_{\Omega_1\times \partial X}\beta_{\xi}^G(X_1(\omega)\inv,(e,n_{X_{2}(\omega)}))\,d\nu_{X_2(\omega)}(\xi)\,d\p_{\overline{m}}(\omega). \\ 
\end{split}
\eeqnn
\eprop 
We use some notation. For $(g,n)\in X$, write $(g,n)\inv:= (g\inv,n)$. The operation\, $\underline{\ast}$ \,for $(g,n), (g',n')\in X$, is $$(g,n)\,\underline{\ast}\,(g',n')=(gg',n).$$ For $\xi\in\partial X$, define $(g,n)\xi=g\xi$.
\bproof
Note that for $i\in\N, i\geq 3$ and $1\leq j\leq i-2$, 
\beqn\label{zerothofintrep}\beta_{\xi}^G(X_i\inv\,\underline{\ast}\,\cdots\,\underline{\ast}\, X_{j+1}\inv\,\underline{\ast}\,X_j\inv,e)=\beta_{X_{j+1}\ast \cdots\ast X_i\, \xi}^G(X_j\inv,(e,n_{X_{j+1}}))\,+\, \beta_\xi^G(X_i\inv\,\underline{\ast}\,\cdots\,\underline{\ast}X_{j+1}\inv,e).\eeqn
We have 
\beqn\label{firstofintrep}
\begin{split}
&\int_{\Omega_1\times\partial X}\frac{\beta_\xi^G(Z_{k}(\omega)\inv,e)}{k}\,d\nu_{(e,n_{X_k(\omega)})}(\xi)\,d\p_{\overline{m}}(\omega) \\ 
&= \int_{\Omega_1\times\partial X}\frac{\beta_\xi^G(X_{k}(\omega)\inv\,\underline{\ast}\,\cdots\,\underline{\ast}\,X_1(\omega)\inv,e)}{k}\,d\nu_{(e,n_{X_k(\omega)})}(\xi)\,d\p_{\overline{m}}(\omega) \\ &
+\int_{\Omega_1\times\partial X} \frac{\beta_\xi^G(Z_{k}(\omega)\inv,X_{k}(\omega)\inv\,\underline{\ast}\,\cdots\,\underline{\ast}\,X_1(\omega)\inv)}{k}\,d\nu_{(e,n_{X_k(\omega)})}(\xi)\,d\p_{\overline{m}}(\omega)
\end{split}
\eeqn
Next for $i\in\N$, set 
\beqnn
\varphi_{i,k}(\omega):=
\begin{cases}
\prod_{j=i+1}^k X_j(\omega)& \text{if}\,i+1\leq  k,\\
e& \text{if}\,i+1> k.\\
\end{cases}
\eeqnn
and note that by the Markov property for $i,k\in\N, i+1\leq k$,
\beqn\label{secondofintrep}
\begin{split}
&\int_{\Omega_1\times\partial X}\frac{\beta_{\varphi_{i,k}(\omega) \xi}^G(X_i(\omega)\inv,(e,n_{X_{i+1}(\omega)}))}{k}\,d\nu_{(e,n_{X_k(\omega)})}(\xi)\,d\p_{\overline{m}}(\omega)\\
&=\frac{1}{k}\sum_{(g_i,n_i)}\sum_{(g_{i+1},n_{i+1})} \sum_{(g,n)}\left[\vphantom{\sum_1}\p_{\overline{m}}[X_i=(g_i,n_i)]\cdot p((e,n_i),(g_{i+1},n_{i+1}))\right.\\
& \times\left.\p_{\overline{m}}[\varphi_{i,k}=(g,n)\,|\,X_{i+1}=(g_{i+1},n_{i+1})]\cdot \int_{\partial X}\beta_{\eta}^G((g_i\inv,n_i),(e,n_{i+1}))\,d\nu_{(g,n)}(\eta)\vphantom{\sum_1}\right]\\
&=\frac{1}{k}\sum_{(g_i,n_i)}\sum_{(g_{i+1},n_{i+1})} \left[ \vphantom{\sum_1} \p_{\overline{m}}[X_i=(g_i,n_i)]\cdot p((e,n_i),(g_{i+1},n_{i+1}))\right. \\
&\times \left. \int_{\partial X}\beta_{\eta}^G((g_i\inv,n_i),(e,n_{i+1}))\,d\nu_{(g_{i+1},n_{i+1})}(\eta)\vphantom{\sum_1}\right] \\
&= \int_{\Omega_1\times\partial X} \frac{\beta_\eta^G(X_i(\omega)\inv,(e,n_{X_{i+1}(\omega)}))}{k}\,d\nu_{X_{i+1}(\omega)}(\eta)\,d\p_{\overline{m}}(\omega).
\end{split}
\eeqn

Let us first consider the first summand on the right of \eqref{firstofintrep}. By \eqref{zerothofintrep} and \eqref{secondofintrep},
\beqn\label{thirdofintrep}
\begin{split}
&\int_{\Omega_1\times\partial X}\frac{\beta_\xi^G(X_{k}(\omega)\inv\,\underline{\ast}\,\cdots\,\underline{\ast}\,X_1(\omega)\inv,e)}{k}\,d\nu_{(e,X_k(\omega))}(\xi)\,d\p_{\overline{m}}(\omega) \\
&= \int_{\Omega_1\times\partial X}\frac{\sum_{i=1}^{k}\beta_{\varphi_{i,k}\xi}^G(X_i\inv,(e,n_{X_{i+1}})) }{k}\,d\nu_{(e,n_{X_k(\omega)})}(\xi)\,d\p_{\overline{m}}\, + \, \int_{\Omega_1\times\partial X}\frac{\beta_\xi^G((e,n_{X_{k+1}}),e)}{k}\,d\nu_{(e,n_{X_k})}(\xi) \,d\p_{\overline{m}}\\
& = \int_{\Omega_1\times\partial X}\frac{\sum_{i=1}^{k}\beta_\xi^G(X_i\inv,(e,n_{X_{i+1}}))}{k}\,d\nu_{X_{i+1}}(\xi)\,d\p_{\overline{m}}\,+\, \int_{\Omega_1\times\partial X}\frac{\beta_\xi^G((e,n_{X_{k+1}}),e)}{k}\,d\nu_{(e,n_{X_k})}(\xi)\,d\p_{\overline{m}},
\end{split} 
\eeqn

The left side of the integral in \eqref{firstofintrep} converges to the drift $l_G$ by hyperbolicity of $X$, harmonicity of $\nu$ and the logarithmic excursion of the walk into cusps. We check the details. First note that $$\int_{\partial X}|\beta_\xi^G(Z_{k}(\omega)\inv,e)|\,d\nu_{(e,n_{X_k(\omega)})}(\xi) \leq \rho_G(e,Z_k).$$  Then by the $L^1$-convergence in Kingman's ergodic theorem
\beqnn
\begin{split}
\int_{\Omega_1}\frac{\rho_G(e,Z_k)}{k}d\p_{\overline{m}} 
&=\int_{\Omega_1}\frac{\rho_G(Z_1,Z_k)}{k}d\p_{\overline{m}}\,+\, O\left(\int_{\Omega_1}\frac{\rho_G(e,Z_1)}{k}d\p_{\overline{m}}\right) \\
&= \int_{\Omega_1}\frac{\rho_G(Z_1,Z_k)}{k}d\p_{\overline{m}}\,+\, O\left(\sum_{j}\frac{j}{k}\cdot \overline{m}((e,j))\right)\\ 
&=  \int_{\Omega_1}\frac{\rho_G(Z_1,Z_k)}{k}d\p_{\overline{m}}\,+ \,O(1/k) \,\longrightarrow\, l_G,
\end{split}
\eeqnn
we have for every $\epsilon>0$ there is $\delta>0$ and $n_\epsilon\in\N$, such that for $A\subset \Omega_1$ with $\p_{\overline{m}}(A)<\delta$ and $k\geq n_\epsilon$,
$$\int_{A\times\partial X}\left|\frac{\beta_\xi^G(Z_{k}(\omega)\inv,e)}{k}\right|\,d\nu_{(e,n_{X_k(\omega)})}(\xi)\,d\p_{\overline{m}}(\omega)<\epsilon.$$ 
Next, by Egorov's theorem, Lemma \ref{tracking} and Lemma \ref{cusp_excur}, there is $\Omega_\epsilon\subset \Omega_1$ measurable, with $\p_{\overline{m}}(\Omega_1\setminus\Omega_\epsilon)<\delta$, such that for all $\omega\in\Omega_\epsilon$ and $k\geq k_\epsilon\in\N$, 
$$\rho_G(Z_1,Z_k)\in ((l_G-\epsilon)\cdot k, (l_G+\epsilon)\cdot k),\quad \rho_X(Z_1,Z_k)\in ((l_X-\epsilon)\cdot k, (l_X+\epsilon)\cdot k), \quad \rho_X(Z_k,\Gamma)\leq C\cdot \log k.$$

For $\omega\in \Omega_1$, by the visuality of $X$, measurably assign a point $\sqrt{Z}_k\inv(\omega)$ on a geodesic ray passing near $Z_k(\omega)\inv$, such that $$\rho_X(e,\sqrt{Z}_k\inv(\omega))=\sqrt{\rho_X(e,Z_k(\omega))}+O(1).$$ Then
\beqnn
\begin{split}
\int_{\partial X}\beta_\xi^G(Z_{k}(\omega)\inv,e)\,d\nu_{(e,n_{X_k(\omega)})}(\xi) 
&=\int_{\partial X\setminus S_X(e,B_X(\sqrt{Z}_k\inv(\omega),R_\Gamma))}\beta_\xi^G(Z_{k}(\omega)\inv,e)\,d\nu_{(e,n_{X_k(\omega)})}(\xi) \\
& + \int_{S_X(e,B_X(\sqrt{Z}_k\inv(\omega),R_\Gamma))}\beta_\xi^G(Z_{k}(\omega)\inv,e)\,d\nu_{(e,n_{X_k(\omega)})}(\xi)\\
\end{split}
\eeqnn
We note for the first integral on the right above that 
\beqnn
\begin{split}
&\int_{\partial X\setminus S_X(e,B_X(\sqrt{Z}_k\inv(\omega),R_X))}\beta_\xi^G(Z_{k}(\omega)\inv,e)\,d\nu_{(e,n_{X_k(\omega)})}(\xi)\\
& = (\rho_X(e,Z_k(\omega))-2\cdot\sqrt{\rho_X(e,Z_k(\omega))}+ O(1))\cdot \nu_{(e,n_{X_k(\omega)})}(\partial X\setminus S_X(e,B_X(\sqrt{Z}_k\inv(\omega),R_X))) 
\end{split}
\eeqnn
For the second integral we note that
\beqnn
\begin{split}
& \int_{S_X(e,B_X(\sqrt{Z}_k\inv(\omega),R_X))}|\beta_\xi^G(Z_{k}(\omega)\inv,e)|\,d\nu_{(e,n_{X_k(\omega)})}(\xi)\\
 &\leq \rho_G(e,Z_k(\omega))\cdot \nu_{(e,n_{X_k(\omega)})}(S_X(e,B_X(\sqrt{Z}_k\inv(\omega),R_X))) \\
& \leq \rho_G(e,Z_k(\omega))\cdot e^{\rho_X(e,(e,n_{X_k(\omega)})}\cdot \nu_{e}(S_X(e,B_X(\sqrt{Z}_k\inv(\omega),R_X)))\\
& = \rho_G(e,Z_k(\omega))\cdot e^{\rho_X(Z_k(\omega),\Gamma)}\cdot \nu_{e}(S_X(e,B_X(\sqrt{Z}_k\inv(\omega),R_X))) \\
&\leq  \rho_G(e,Z_k(\omega))\cdot e^{\rho_X(Z_k(\omega),\Gamma)}\cdot e^{-\sqrt{\rho_X(e,Z_k(\omega))}}=o(k).
\end{split}
\eeqnn
It follows that 
$$\lim_{k\to\infty}\int_{\Omega_\epsilon\times\partial X}\frac{\beta_\xi^G(Z_{k}(\omega)\inv,e)}{k}\,d\nu_{(e,n_{X_k(\omega)})}(\xi)\,d\p_{\overline{m}}=l_G\cdot \p_{\overline{m}}(\Omega_\epsilon).$$
The second summand in the right of \eqref{firstofintrep} and the second summand in \eqref{thirdofintrep} similarly converge to zero  as the depth is at most logarithmic in time. 

The first summand in the right side of the last equality of \eqref{thirdofintrep} converges to the right side of the integral in the claim of the proposition by the pointwise ergodic theorem ($(\Omega_1,T,\p_{\overline{m}})$ is ergodic). The proposition follows.
\eproof

\section{Dimensions of $\nu$ and $\mu$}
We will use the following corollary of Lemma \ref{cusp_excur} (cf. \cite[Proposition 4.10]{SV}) which gives uniform estimates for the Patterson-Sullivan measure of sufficiently small balls centered in sets from an increasing sequence of compact sets which exhaust the boundary. As the harmonic measure is Ahlfors-regular, a similar but stronger conclusion for it is obvious.
\bcor\label{exhaust}
There are constants $C\geq 1$ and $\alpha_1, \alpha_2, R>0$, such that for $x\in\Gamma\subset X$, there exists a nested sequence of compact sets $E_i\subset\partial X$, for which $\partial X=\bigcup_{n=1}^{\infty} E_n$ (mod $\mu$), and for all $n\in\N$, and for all $\xi\in E_n$, there exist $t_n>0$, such that for all $t>t_n$ (where $\gamma^x(t)\in X$), if $\gamma^x$ is a geodesic starting at $x$, with $\gamma^x_\infty=\xi$, the following holds: 
$$C_1\cdot t^{-\alpha_1\cdot (1+1/n)}\cdot e^{-D_\Gamma\cdot t}\leq \mu_x(S(x,B(\gamma^x(t), R)))\leq C_2\cdot t^{\alpha_2\cdot (1+1/n)}\cdot e^{-D_\Gamma\cdot t}.$$
\ecor

\bproof
Apply the measurable selection theorem to associate measurably (with respect to $\mu$) to each $\xi\in\partial X$, a geodesic ray $\gamma^{x,\xi}$ starting at $x$, such that $\gamma^{x,\xi}(\infty)=\xi$.
Given $n\in\N$, by Luzin's theorem applied to the function $$t_n(\xi)=\inf\,\left\lbrace t>0\,\mid\, \text{for}\;t'>t,\;\rho_X(\gamma^{x,\xi}(t'),\Gamma)\leq \log t'\cdot \frac{1+1/n}{2\cdot D_\Gamma-d_H\cdot \log a}\right\rbrace,$$ (which is finite $\nu_x$-a.e $\xi\in\partial X$ by Lemma \ref{cusp_excur}), 
obtain a compact set $E_n$ such that $$\nu_X(\partial X\setminus E_n)<1/n,$$ and where the function $t_n$ is bounded by a number also written $t_n$.
Then for $\xi\in E_n$, for $t>t_n$, if $\gamma^{x,\xi}(t)\in \Gamma$, the claim holds. So assume that $\gamma^{x,\xi}(t)\in X\setminus\Gamma$. 
Then, choosing $R$ large enough, we have by Remark \ref{imp_1},
\beqnn
\mu_x(S(x,B_X(\gamma^{x,\xi}(t),R)))\approx e^{-D_\Gamma\cdot t}\cdot e^{-D_\Gamma\cdot\rho_X(\gamma^{x,\xi}(t),\Gamma)\cdot \left(1-\frac{k_H(\gamma^{x,\xi}(t))}{D_\Gamma\cdot\rho_X(\gamma^{x,\xi}(t),\Gamma)}\right)},
\eeqnn
and hence since $\xi\in E_n$,
\beqnn
e^{-D_\Gamma\cdot t}\cdot \left(t^{-\frac{1+1/n}{2-d_H\cdot\log a/D_\Gamma}}\right)^{1+d_H\cdot\log a/D_\Gamma}\lesssim\; \mu_x(S(x,B_X(\gamma^{x,\xi}(t),R)))\, \lesssim  e^{-D_\Gamma\cdot t}\cdot \left(t^{\frac{1+1/n}{2-d_H\cdot\log a/D_\Gamma}}\right)^{1+d_H\cdot\log a/D_\Gamma}.
\eeqnn
The lemma follows since any other $(X,\rho_X)$-geodesic ray $\gamma^x$ with end point $\xi$ is $\delta_X$-close to $\gamma^{x,\xi}$.
\eproof

\bproof[\textbf{Proof of Theorem \ref{dimdriftent}}]: 
Consider a set $\Omega'\subset X^{\N}$ of full measure using Lemma \ref{drift} and Lemma \ref{tracking}, where the following hold for all $\omega\in\Omega'$, $k>n_\omega$:
\beqn\label{firstofdd}
\left|\frac{\rho_G(Y_0(\omega),Y_k(\omega))}{k}-l_G\right|\leq \eta,\quad\left|\frac{\rho_X(Y_0(\omega),Y_k(\omega))}{k}-l\right|\leq \eta,
\eeqn and 
\beqn\label{secondofdd}
\rho_G(Y_k(\omega),\gamma^{x,Y_\infty(\omega)})\leq 3\cdot \log k.
\eeqn
Write $x_n(\omega)$ for measurably assigned points for $\omega\in \Omega'$ and $n\in\N$, which are nearest projections from $Y_n(\omega)$ to $\gamma^{x,Y_\infty(\omega)}$. From \eqref{firstofdd} and \eqref{secondofdd}, we have
$$
|\rho_G(Y_0(\omega),x_k(\omega))-l_G\cdot k|\leq 2\cdot\eta\cdot k,\quad \text{for all}\;k>n'_\omega.
$$ 
For $\p_e$-a.e $\omega\in\Omega'$, $\xi_\omega=Y_\infty(\omega)$ and $k>n'_\omega$, note that by Lemma \ref{exhaust}, setting $$r_k=e^{-\epsilon_X\cdot\rho_X(e,x_k(\omega))},$$ we have
$$
e^{-(l_G + 2\cdot\eta)\cdot k}\lesssim\nu_e(S_G(e, B_G(x_k(\omega),R_\nu)))\lesssim  e^{-(l_G - 2\cdot\eta)\cdot k},
$$
and thus 
$$
e^{-(l_G + 2\cdot\eta)\cdot k}\lesssim \nu_e(B_\infty^X(\xi_\omega, r_k))\lesssim e^{-(l_G - 2\cdot\eta)\cdot k}, 
$$ from which we get, by letting $\eta\to 0$, for $\omega\in\Omega'_n$,
$$\lim_{k\to \infty}\frac{\log \nu_e(B_\infty^X(\xi_\omega, r_k))}{\log r_k}=\frac{l_G}{\epsilon_X\cdot l}.$$ The theorem follows for the base-point $e\in\Gamma$ by letting $n\to\infty$, using the doubling property of $\nu$ (see \cite[p. 26]{BHM}). For general basepoint, one uses the harmonicity of $\nu$. The existence of the second limit follows similarly. 

For $\mu_x$, $x\in\Gamma$, consider for $n\in\N$, the set $E_n$. For $\xi\in E_n$, we get$$\lim_{r\to\infty}\frac{\log\mu_x(B_{x,\infty}(\xi,r))}{\log r}=-\lim_{n\to\infty}\frac{\log \mu_x(S_X(x,B_X(\gamma^{x,\xi}(n),r)))}{\epsilon_X\cdot n}=\frac{D_\Gamma}{\epsilon_X};$$ and $E_n$ exhaust $\partial X$ $(\text{mod}\,\mu)$. This finishes the proof of the first and third parts of the claim. 

The approach of \cite{LeP} and \cite{Tan} can be adapted to our setting to obtain the proof of the second part the claim. The key ingredients towards that end are Theorem \ref{drift}, hyperbolicity of $X$, ergodicity of $(X,\Gamma,\nu)$ and the following (recall $(\cdot|\cdot)_{\cdot,G}$ and $(\cdot|\cdot)_{\cdot,X}$ are the Gromov products in $\rho_G$ and $\rho_X$ respectively).
\blem We have for $x\in X$,
$$
\lim_{n\to\infty}\frac{(Y_n(\omega)|Y_{n+1}(\omega))_{x,G}}{n}=l_G,\quad \lim_{n\to\infty}\frac{(Y_n(\omega)|Y_{n+1}(\omega))_{x,X}}{n}=l,
$$
$\p_x$-a.e $\omega\in X^{\N}$.
\elem
\bproof
This is a consequence of Theorem \ref{drift} and Lemma \ref{tracking}.
\eproof
The proof of the second part (dimension in terms of entropy) of Theorem \ref{dimdriftent} then follows closely the steps in \cite{LeP} and \cite{Tan}. The upper bound on the dimension uses the same argument. The lower bound needs slight modification using the asymptotic behaviour or $m(Y_n)$ as in Lemma \ref{fund}. We omit repeating the arguments.
\eproof

\vspace{0.5cm}
Random walks with suitable classes of infinite range transition probabilities, comparison of the harmonic measure with some natural classes of measures on the Bowditch boundary in this case, random walks on weak cusped graphs and consequences are topics of upcoming work.

\end{document}